\documentclass[aoas,preprint]{imsart}

\RequirePackage[OT1]{fontenc}
\RequirePackage{amsthm,amsmath,amssymb,enumerate}
\RequirePackage[authoryear]{natbib}
\RequirePackage[colorlinks,citecolor=blue,urlcolor=blue]{hyperref}
\usepackage{algorithm}
\usepackage{algcompatible}
\usepackage{caption}
\usepackage{subcaption}
\usepackage{graphicx}
\usepackage{geometry}
\startlocaldefs
\numberwithin{equation}{section}
\theoremstyle{plain}
\newtheorem{thm}{Theorem}
\newtheorem{lemma}{Lemma}
\newtheorem{corollary}{Corollary}				
\newtheorem{proposition}{Proposition}
\theoremstyle{definition}
\newtheorem{definition}{Definition}
\newtheorem{assumption}{Assumption}
\newtheorem{remark}{Remark}

\def\P{{\mathbb P}}
\def\R{{\mathbb R}}

\def\N{{\mathbb N}}

\def\cF {\mathcal{F}}

\endlocaldefs

\begin{document}

\begin{frontmatter}
\title{Stability,  convergence to equilibrium and simulation of non-linear Hawkes Processes with memory kernels given by the sum of Erlang kernels}
\runtitle{}

\begin{aug}
\author{\fnms{A.} \snm{Duarte}\thanksref{m2}\ead[label=e1]{alineduarte@usp.br}},
\author{\fnms{E.} \snm{L\"ocherbach}\thanksref{m1}\ead[label=e2]{eva.loecherbach@u-cergy.fr}}
\and
\author{\fnms{G.} \snm{Ost}\thanksref{m3}
\ead[label=e3]{guilhermeost@im.ufrj.br}}

\runauthor{A. Duarte et al.}

\affiliation{Universit\'e de Cergy-Pontoise\thanksmark{m1}, 
Universidade de S\~ao Paulo\thanksmark{m2} and Universidade Federal do Rio de Janeiro \thanksmark{m3} }

\address{A. Duarte\\
Universidade de S\~ao Paulo,\\ 
Instituto de Matem\'atica e Estat\'istica,\\
S\~ao Paulo,\\
Brazil.\\
\printead{e1}}

\address{E. L\"ocherbach\\
Universit\'e de Cergy-Pontoise,\\ 
AGM, CNRS-UMR 8088,\\
95000 Cergy-Pontoise,\\
France.\\
\printead{e2}}

\address{G. Ost\\
Universidade de Federal do Rio de Janeiro,\\ 
Instituto de Matem\'atica, \\
Rio de Janeiro,\\
Brazil.
\\
\printead{e3}}

\end{aug}

\begin{abstract}
Non-linear Hawkes processes with memory kernels given by the sum of Erlang kernels are considered. 
It is shown that their stability properties can be studied in terms of an associated class of piecewise deterministic Markov processes, called {\it Markovian cascades of successive memory terms}. Explicit conditions implying the positive Harris recurrence of these processes are presented. The proof is based on integration by parts with respect to the jump times. A crucial property is the non-degeneracy of the transition semigroup which is obtained thanks to the invertibility of an associated Vandermonde matrix. For Lipschitz continuous rate functions we also show that these Markovian cascades converge to equilibrium exponentially fast with respect to the Wasserstein distance. Finally, an extension of the classical thinning algorithm is proposed to simulate such Markovian cascades. 
\end{abstract}



\end{frontmatter}

\section{Introduction}
Hawkes processes have regained a lot of interest in the recent years, in particular in econometrics, as good models to account for contagion risk and clustering arrival of events. They have shown to be very useful also in neuroscience due to their capacity of reproducing the typical time dependencies observed in spike trains of neurons as well as the interaction structure of neural nets. Originally introduced by  \cite{Hawkes1971} and \cite{ho} as a model for the appearances of earthquakes, their key feature is the fact that any event is able to trigger future events -- for this reason, Hawkes processes are sometimes called ``self-exciting point processes''. In their by now classical paper, \cite{bm} develop the stability theory of general non-linear Hawkes processes, also in a multivariate frame. \cite{hrbr} have put the foundations for the use of Hawkes processes as models of spike trains in neuroscience, see also \cite{ccdr}, and recently some effort has been spent to study Hawkes processes in high dimensions, especially focusing on properties such as the propagation of chaos, see \cite{dfh} and  \cite{chevallier}, see also \cite{SusanneEva} in a multi-class frame. Finally, we refer to \cite{zhu2015} for a study of the large deviation properties of non-linear Hawkes processes having Markovian intensity function.

In the present paper Hawkes processes with memory kernels given by the sum of Erlang kernels are considered. It is shown that the longtime behaviour and stability properties of these processes can be studied in terms of an associated class of piecewise deterministic Markov processes (PDMPs). More precisely, let $N$ be a counting process on $\R_+$ characterised by its intensity process $(\lambda_t)_{t\geq 0}$ defined, for each $t\geq 0$,  through the relation 
$$ \P ( N \mbox{ has a jump in } ]t , t + dt ] | {\cal F}_t ) = \lambda_t  dt , $$
where $ {\cal F}_t = \sigma (  N ( ]u, s ] ) , \, 0\leq u <  s \le t )$ and 
\begin{equation}\label{eq:intensity0}
\lambda_ t = f  \left(\delta+ \int_{]0  , t [} h ( t-s) d N_s  \right) .
\end{equation}
Here, $f : \R \to \R_+$ is the {\it jump rate function} and $ h : \R_+ \to \R  $ is the {\it memory kernel.} 
The parameter $\delta\in\R$ is interpreted as an initial input to the jump rate function.

The memory kernel $h$ is assumed to be given by the sum of Erlang kernels, that is, $ h = \sum_{i=1}^L h_i, $ where each function $h_i $ is of the form 
\begin{equation}\label{eq:Erlang}
h_i(t)=c_ie^{-\alpha_i t}\frac{t^{n_i}}{n_i!}, t \geq 0,
\end{equation}
where $ c_i \in \R , \alpha_i > 0 $ and $n_i  \in \N .$ 

Erlang kernels are widely used in the modelling literature, for example to model delays in the hemodynamics in nephrons, see \cite{Ditlevsen2005,KidneyDelay} or to prove the existence of oscillations in large-scale limits of interacting neurons in a mean-field frame, see \cite{SusanneEva}. This is the main motivation for the particular choice taken for the memory kernel $h.$ Moreover, it is well known that the class of memory kernels having this form is dense in $ L^1 ( \R_+),$
see e.g.\ \cite{kammler}. Therefore, any Hawkes process having integrable memory kernel can be well approximated by a Hawkes process having an Erlang memory kernel, see \eqref{eq:approx} below, at least over compact time intervals. Finally,  the specific memory structure induced by Erlang kernels allows a completely new approach to simulation of non-linear Hawkes processes.

Erlang kernels depend on three parameters, $c_i, n_i $ and $\alpha_i .$ Here, $n_i+1$ is the {\it order} of the delay of the influence of a past event on a future event. It takes its maximum absolute value at $(n_i+1)/\alpha_i$ time units back in time. The mean is $(n_i+1)/\alpha_i$ (if normalising to a probability density). The higher the order of the delay, the more concentrated is the delay around its mean value, and in the limit of $n_i \rightarrow \infty$ while keeping $(n_i+1)/\alpha_i$ fixed, the delay converges to a discrete delay. The sign of $c_i$ indicates if the influence of past events on future events is inhibitory or excitatory.

The use of Erlang kernels allows to relate the study of the longtime behaviour of a Hawkes process having intensity \eqref{eq:intensity0} to the study of an associated system of PDMPs. 
More specifically, it is easily shown that the system of stochastic processes $X^{(i,0) }_t= \int_{]0, t]} h_i  ( t-s) d N_s, t\geq 0$, for each $1 \le i \le L, $ can be completed, by introducing $ \sum_{i=1}^L n_i $ auxiliary processes, to a piecewise deterministic Markov process in dimension $ L + \sum_{i=1}^L n_i  , $ see \eqref{eq:cascade} below. Between successive jumps of $N, $ the evolution of each $X^{( i,0 )}_t,$ together with its auxiliary processes, is explicitly given by a deterministic flow. Jumps do only occur in the auxiliary variables $ X^{(i, n_i ) }, 1 \le i \le L . $ We shall call this class of PDMPs {\it Markovian cascades of successive memory terms}.

We prove that these Markovian cascades are recurrent in the sense of Harris under the usual sub-criticality condition
\begin{equation}
\label{cond:Bremaud}
\| f \|_{Lip} \int_0^\infty | h ( t)| dt < 1. \footnote{if all $ \alpha_i $ are equal and all $ c_i$ are of the same sign. In the case of bounded rate functions, we do not need to impose this condition.} 
\end{equation} 

Under \eqref{cond:Bremaud}, we are able to construct a Lyapunov function implying that these processes return to a compact set infinitely often, almost surely. Under the additional condition of some minimal {\it ellipticity}, that is, some minimal jump activity, we establish, in Theorem \ref{thm:Doeblin}, a Doeblin like lower bound based on integration by parts with respect to the jump times. A crucial property is the non-degeneracy of the transition semigroup which is obtained thanks to the invertibility of an associated Vandermonde matrix and structure of the flow of the Markovian cascades (see \eqref{eq:flow} below). In the case of Lipschitz continuous rate functions we also show that the Markovian cascades converge to equilibrium exponentially fast with respect to the Wasserstein distance. 

The fact that the flow governing the evolution of the Markovian cascades in between successive jump times is explicitly given enables us to introduce an efficient simulation algorithm which allows to sample from $N $ on $  [0, T ] $ for any finite time horizon $T > 0 $ and any fixed parameter $\delta\in \R$. This method is straightforward to implement, and can be easily extended to multi-dimensional versions.  

The paper is organised as follows. In Section \ref{sec:def}, we present the model and provide some preliminary remarks. In Section \ref{sec:long_time_term}, the long-time behaviour of the Markovian cascades is investigated. 
The statement and proof for the case $L=1$ of Theorem \ref{thm:Doeblin}, establishing the Doeblin lower bound for the Markovian cascades, is also included in this section. 
In Section \ref{sec:Simulations}, a simulation algorithm to simulate simultaneously a Hawkes process with memory kernel given by the sum of Erlang kernels and its Markovian cascade is proposed. In Section \ref{sec:numerical_examples}, numerical examples are presented. Finally, in the Appendix \ref{Appendix}, we prove Theorem \ref{thm:Doeblin} in the general case.
 
\section{Model definition and preliminary remarks}\label{sec:def}
Throughout the article the set $\N$ denotes the set of non-negative integers, $\N^*$ the set of positive integers $\{1,2,\ldots\}$ and 
$\mathcal{B}((0,\infty))$ (resp.
$\mathcal{B}((a,b])$, for real numbers $0\leq a\leq b<\infty$) the  Borel sigma-algebra on $(0,\infty)$ (resp. on $(a,b]$). 

We work on the following filtered space $(\Omega,\cF, \mathbb{F})$. Let $\Omega$ be the canonical path space of simple point processes, i.e., 
$$
\Omega=\{w=(t_n)_{n\in \N^*}\in ] 0, \infty ]^{\N^*}:   t_n\leq t_{n+ 1}, t_n<t_{n+ 1}\ \mbox{if} \ t_n<+\infty, \ \lim_{n\to +\infty} t_n =+\infty\}.
$$
For each $w\in \Omega$ and $n\in\N^*$, we define $T_n(w)=t_n$. For each $w\in \Omega$, we associate the canonical point measure  $\mathcal{B}((0,\infty)) \ni A \mapsto N(w)(A)=\sum_{n\in \N^*}\delta_{T_n(w)}(A)$. We shall write for short $N(A)$ rather than $N(w)(A)$; when $A=(0,t]$ for some $t\geq 0$,  we simply write $N_t$ to denote $N((0,t])$. Finally, we define $\cF_t=\sigma\{ N (A):A\in \mathcal{B}((0,t])\}$ for each $t\geq 0$, $\cF=\sigma\{ N (A):A\in \mathcal{B}((0,\infty))\}$ and $\mathbb{F}=(\cF_t)_{t\geq 0}$. 

Let $ f : \R \to \R_+$ and $ h : \R _+ \to \R$ be measurable functions and let $ n$ be a deterministic point process on $ ]- \infty, 0 ] $ such that $ \int_{ ] - \infty , 0]} h (t- s) n (ds) $ is finite for all $t \geq 0.$ 
\begin{definition}
A Hawkes process with parameters $(f,h)$ and with initial condition $n$ is a probability measure $P$ on the filtered space $(\Omega,\cF, \mathbb{F})$ such that the compensator of $(N_t)_{t\geq 0} $  is given by
$
(\int_0^t\lambda_s ds)_{t\geq 0}, 
$
where $(\lambda_t)_{t\geq 0}$ is the non-negative $\mathbb{F}-$predictable process defined for $t\geq 0$ by
\begin{equation}\label{eq:intensity}
\lambda_t=f\left( \int_{ ] - \infty , 0]} h (t- s)  n (ds) +\int_{]0, t[ }h(t-s)dN_s\right).
\end{equation} 
\end{definition}

The stochastic process $(\lambda_t)_{t\geq 0}$ is called {\it intensity process}. The functions $f : \R \to \R_+ $ and $h: \R_+ \to \R$ are called {\it jump rate function} and {\it memory kernel} respectively. We shall work under the following assumptions.

\begin{assumption}\label{ass:Lip}
The rate function $f : \R \to \R_+$ is either bounded or Lipschitz continuous with Lipschitz constant $ \|f\|_{Lip}.$
\end{assumption}

\begin{assumption}\label{ass:h}
The memory kernel $h: \R_+ \to \R$ can be written as sum of Erlang kernels, i.e, for each $t\geq 0$,
\begin{equation}\label{eq:erlang}
h(t)=\sum_{i=1}^L c_ie^{-\alpha_i t}\frac{t^{n_i}}{n_i!},
\end{equation}
where for each $1\leq i\leq L$, $c_i\in \R$, $\alpha_i>0$ and $n_i\in \N$.
\end{assumption}

 

Under Assumption \ref{ass:h}, the intensity process \eqref{eq:intensity} can be described by an associated PDMP. Indeed, for each $ 1 \le i \le L$ and  $0\leq k\leq n_i$, writing  for each $t\geq 0$, 
\begin{equation}\label{eq:Xik}
X^{(i, k)}_t= \int_{ ]- \infty , 0] } c_ie^{-\alpha_i (t-s)}\frac{(t-s)^{(n_i-k)}}{{(n_i-k)}!} n (ds) +\int_{]0 ,t]}c_ie^{-\alpha_i (t-s)}\frac{(t-s)^{(n_i-k)}}{{(n_i-k)}!}dN_s,
\end{equation}
we have 
$$ \lambda_t = f\left(X^{(1,0)}_{t-} + \ldots + X^{ ( L, 0) }_{t-}\right) ,$$
and one easily deduces that for $t\geq 0$ and $ 1 \le i \le L,$
\begin{eqnarray}\label{eq:cascade}
dX^{(i, 0)}_t&=&X^{(i, 1)}_tdt-\alpha_i X^{(i, 0)}_tdt \\
 &\vdots& \nonumber \\
dX^{(i, n_i-1)}_t&=&X^{(i,n_i)}_tdt-\alpha_i X^{(i,n_i-1)}_tdt  \nonumber\\
dX^{(i, n_i)}_t&=&-\alpha_i X^{(i, n_i)}_tdt+c_i dN_t,  \nonumber
\end{eqnarray}
with initial condition $X^{(i,k)}_0=x_0^{(i,k)} =  \int_{ ]- \infty , 0] } c_ie^{\alpha_i s}\frac{(-s)^{(n_i-k)}}{{(n_i-k)}!} n (ds)$. 

\begin{figure}[h!]
\includegraphics[scale=0.65]{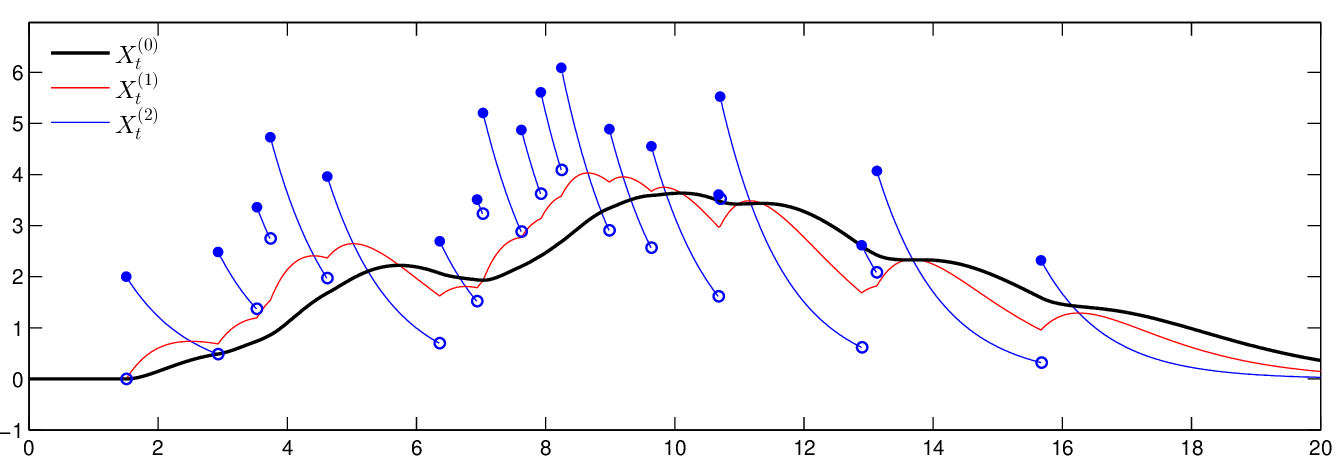}
\includegraphics[scale=0.65]{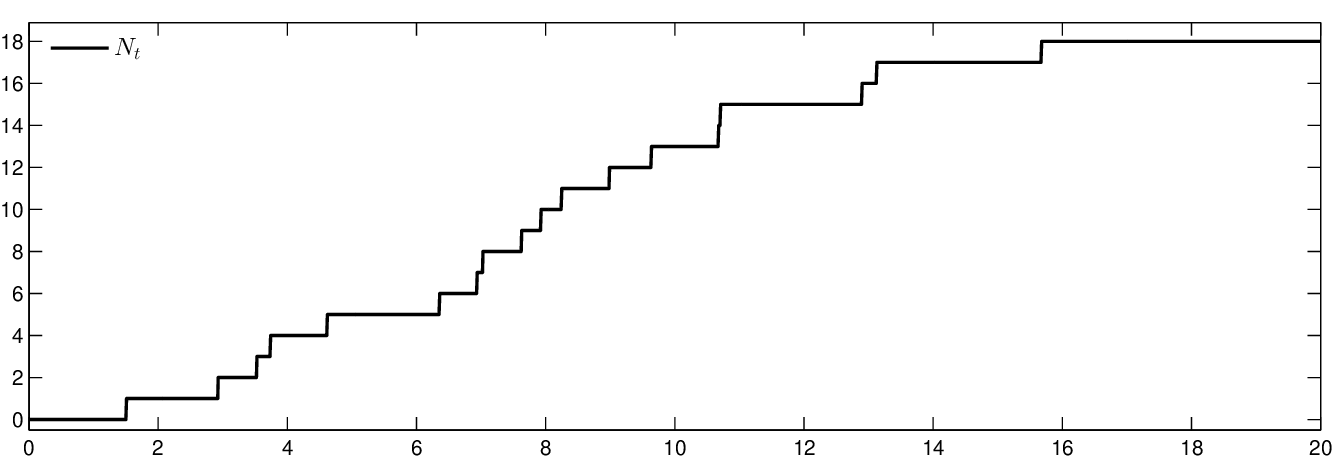} 
\caption{A finite joint realisation of the Markovian cascade $X=(X_t)_{0\leq t\leq T}$ (upper panel) and its associated counting process $N=(N_t)_{0\leq t\leq T }$ (lower panel) for the choices $L=1,$ $n_1=2$, $c_1=2$, $\alpha_1=1$, $T=20$ and $f(x)=x/5+1$ with initial input $x_0=(x^{(0)}_0,x^{(1)}_0,x^{(2)}_0)=(0,0,0)$. The blue (resp. red and black) trajectory corresponds to the realisation of $(X^{(2)}_t)_{0\leq t\leq T}$ (resp. $(X^{(1)}_t)_{0\leq t\leq T}$ and $(X^{(0)}_t)_{0\leq t\leq T}$). Notice that the smaller the index $0\leq k\leq n_1 $ the smoother the correspondent process $(X^{(k)}_{0\leq t\leq T}).$ }
\label{sample_figure}
\end{figure}

Write $ \kappa = L + \sum_{ i=1}^L n_i .$ The associated PDMP is the Markov process $ X=(X_t)_{t\geq 0}$ having c\`adl\`ag paths and taking values in $\R^{\kappa },$ defined, for each $t\geq 0$, by 
\begin{equation}
\label{def:PDMP}
X_t=\Big(X^{(1)}_t,\ldots, X^{(L)}_t\Big) \ \mbox{with} \ X^{(i)}_t=\Big(X^{(i,0)}_t,\ldots, X^{(i,n_i)}_t\Big), \ 1\leq i\leq L.
\end{equation}
If $L=1,$ that is, $h$ is a pure Erlang kernel, we write for short $ X_t = ( X^{ (k)}_t , 0 \le k \le n_1).$
We call the process $X$ {\it Markovian
cascade of successive memory terms}. Its infinitesimal generator $\mathcal{L}$ is given for any smooth test function $g:\R^{\kappa}\mapsto \R$ by
\begin{equation}\label{eq:gen}
\mathcal{L}g(x)=\langle F(x),\nabla g(x)\rangle +f\Big(\sum_{i=1}^L x^{(i,0)}\Big)\Big(g\big(x+\sum_{i=1}^L c_i e_{(i,n_i)}\big)-g(x)\Big),
\end{equation}
where $x=\big(x^{(1)},\ldots, x^{(L)}\big)$ with $x^{(i)}=(x^{(i,0)},\ldots, x^{(i,n_i)})$ and
$e_{(i,n_i) }\in \R^{\kappa }$ is the unit vector having entry $ 1 $ in the coordinate $ (i, n_i), $ and $0$ elsewhere. Finally, $F:\R^{\kappa}\mapsto \R^{\kappa }$ is the vector field associated to the system of first-order ODE's
\begin{equation}
\label{def:sysODEs}
\begin{cases} 
\dfrac{d}{dt}x^{(i, 0)}_t=x^{(i,1)}_t-\alpha_i x^{(i,0)}_t \\
\vdots \\
\dfrac{d}{dt}x^{(i, n_i-1)}_t=x^{(i, n_i)}_t-\alpha_i x^{(i,n_i-1)}_t 
\vspace{0.15cm}\\
\dfrac{d}{dt}x^{(i, n_i)}_t=-\alpha x^{(i, n_i)}_t ,  \ 1 \le i \le L, 
\end{cases}
\end{equation}
given by $F(x) = ( (F^{(1)} (x), \ldots , F^{(L)} ( x) ) ,$ where $ F^{(i)} (x) =  (F^{(i,0)} (x) , \ldots , F^{(i,n_i)} (x) )$ 
with
$$ F^{(i,k)} (x) = - \alpha_i x^{(i,k)} + x^{(i,k+1)}\ \mbox{for} \ 0 \le k < n_i, \ \mbox{and} \ F^{(i,n_i)} (x) = - \alpha_i x^{(i,n_i)}.$$
Notice that jumps introduce discontinuities only in the coordinates $X^{(i,n_i)}_t$ of $X_t$. Figure \ref{sample_figure} depicts a realisation of the joint processes $(N_t, X_t)_{t\geq 0}$ in the case $L=1.$ 

Hereafter, we write $\varphi_t(x)=(\varphi^{(1)}_t(x),\ldots, \varphi^{(L)}_t(x))$
for the unique solution, starting from $x\in \R^{\kappa }$, of the system \eqref{def:sysODEs}. It is immediate to check that for each $1 \le i \le L$ and $t\geq 0$,
\begin{equation}\label{eq:flow}
\varphi^{(i)}_t(x)=
\left(
\begin{array}{c}
e^{-\alpha_i t}(x^{(i,0)}+tx^{(i,1)}+\ldots+ \frac{t^{n_i}}{n_i!}x^{(i,n_i)})\\
\vdots\\
e^{-\alpha_i t}(x^{(i,n_i-1)}+tx^{(i,n_i)})\\
e^{-\alpha_i t}x^{(i,n_i)}
\end{array}
\right) .
\end{equation}
Notice that $\varphi^{(i)}_t(x)$ depends only on the variable $x^{(i)}$.
Given the Markovian cascade of successive memory terms \eqref{eq:cascade}-\eqref{eq:gen}, one recovers immediately the non-linear Hawkes processes with intensity \eqref{eq:intensity} as shows the following proposition.
In what follows, for any $x_0\in\R^{\kappa}$, we write $P_{x_0}$ for the probability measure on $\R^{\kappa}$ under which $X_0=x_0 ,$ and denote by $E_{x_0}$ the expectation taken with respect to $P_{x_0}$.  


\begin{proposition}
\label{thm:1}
Suppose Assumptions \ref{ass:Lip} and \ref{ass:h}. Fix any initial condition $n$ on $  ] -  \infty, 0 ] $ such that $ \int_{ ]- \infty , 0 ]} h(t-s )  n (ds)  < \infty $ for all $ t \geq 0$ and put, for each $1\leq i\leq L$ and $0\leq k\leq n_i$, $x^{(i,k)}_0 = \int_{ ]- \infty , 0] } c_ie^{\alpha_i s}\frac{(-s)^{(n_i-k)}}{{(n_i-k)}!} n (ds) . $ 
Let $X= (X_t)_{t \geq 0}$ be the Markov process having generator \eqref{eq:gen}, starting from $X_0=x_0.$ Then $X$ is non-explosive, i.e., $X$ has $P_{x_0}$-almost surely a finite number of jumps on each interval $ [s, t ], 0 \le s < t < \infty .$ 
Finally, introduce $ N_t = \sum_{ s \le t } 1_{\{\Delta X_s \neq 0\}}$ the counting process associated to the jumps of $X.$ Then $(N_t)_{t\geq 0}$ is a non-linear Hawkes process with intensity \eqref{eq:intensity}. 
\end{proposition}

\begin{proof}
Let us define $\mu (dt)= \sum_{ n \geq 1 } \delta_{T_n } (d t) $ and $\nu (dt ) = f \big(\sum_{i=1}^L X^{(i, 0)}_t\big) dt .$ Then Proposition 3.1 in \cite{Jacod75} implies that $\nu$ is the predictable compensator of $ \mu .$ In particular, the compensator of $N_t$ is given by $\int_0^t \lambda_s ds $ with $ \lambda_s = f \big(\sum_{i=1}^L X^{(i, 0)}_{s-}\big) .$ 

It remains to prove the non-explosiveness of the process $X.$ In the case of bounded $f$, nothing has to be proved. Suppose therefore that $f$ is Lipschitz continuous and define $g(x)= f(0) (\sum_{i=1}^L |c_i|) +\sum_{i=1}^L \sum_{k=0}^{n_i} |x^{(i,k)}|$ for $x\in\R^{\kappa }.$ Let $ A= \max_{1\leq i\leq L} \alpha_i $ and $ c= \sum_{i=1}^L |c_i |.$ 
By plugging $g$ in \eqref{eq:gen}, we have
\begin{eqnarray*}
Lg(x)&\leq & \sum_{i=1}^L\Big( \sum_{k=0}^{n_i-1} sg(x^{(i,k)})\big(x^{(i, k+1)}-\alpha_i x^{(i, k)}\big)-\alpha_i |x^{(i, n_i)}|+|c_i|f\big(\sum_{j=1}^L x^{(j, 0)}\big)\Big)\\
&\leq &c f(0)+( c \; \|f\|_{Lip} + \alpha) \sum_{i=1}^L | x^{(i, 0)}|+(1+ A)\sum_{ i=1}^L \sum_{k=1}^{n}|x^{(i, k)}|\\
&\leq & Cg(x),
\end{eqnarray*}
where $C= C (A , c ,  \|f\|_{Lip} )$ and $sg(y)$ is the sign of $y\in \R$. In the second inequality above we have used that $f(y)\leq \|f\|_{Lip}|y|+f(0)$ for any  $y\in \R$. 
Thus, by applying Dynkin's formula  and then using $Lg(x)\leq C g(x)$, one concludes that
$$
E_{x_0}[g(X_t)]=g(x_0)+\int_{0}^tE_{x_0}[Lg(X_s)]ds
\leq  g(x_0)+C\int_{0}^tE_{x_0}[g(X_s)]ds.
$$
Then by Gronwall's inequality, $E_{x_0}[g(X_t)]\leq g(x_0)e^{Ct}$. From this last estimate we conclude the proof noticing that
$$ 
E_{x_0}( N( [s, t]) ) = E_{x_0}\int_s^t f \big(\sum_{i=1}^L X^{(i, 0)}_u\big)  du \le \| f \|_{Lip} \int_s^t  \sum_{i=1}^L E_{x_0} | X^{(i, 0)}_u | du + f(0 ) (t-s).
$$
Since $\sum_{i=1}^L | X^{(i, 0)}_u | \le g ( X_u )$, it follows immediately from the inequality above that
\begin{eqnarray*}
E_{x_0} ( N( [s, t]) )&\leq &\|f\|_{Lip}\int_s^t E_{x_0}[g(X_u)]du+f(0)(t-s)\\
&\leq & g(x_0) \frac{\|f\|_{Lip}}{C}(e^{Ct}-e^{Cs})+f(0)(t-s)<\infty .
\end{eqnarray*} 
\end{proof}

\begin{remark}
The converse statement of the above proposition does also hold true. More precisely, let $(N_t)_{t\geq 0}$ be a non-linear Hawkes process with intensity \eqref{eq:intensity} where $h $ is given by \eqref{eq:erlang}. Suppose moreover that $N$ starts from the  $ n (ds) $ on $\R_{-}, $ where $n$ is some discrete point measure on $\R_{-}$ such that $ x^{(i, k)}_0 := \int_{-\infty}^0 c_ie^{-\alpha_i (-s)}\frac{(-s)^{(n_i-k)}}{{(n_i-k)}!} n (ds) $ are well-defined. Introduce the associated processes $ X_t^{(i,k ) }, 1 \le i \le L,  0  \le k \le n_i $ as in \eqref{eq:Xik}. Then $X= (X_t)_{t \geq 0}$ is Markov with generator \eqref{eq:gen}. 
\end{remark}

\subsection{Some comments on the use of Erlang kernels}
Erlang kernels are widely used in the modeling literature. They have been introduced by Erlang in the 1920's to provide an efficient approach for analyzing telephone networks.  Nowadays, they are widely used in the theoretical and mathematical biology literature, see e.g. \cite{Ditlevsen2005,KidneyDelay} where they serve as a good model to describe delays in the hemodynamics in nephrons. They are also the building block to prove the existence of oscillations in large-scale limits of interacting neurons in a mean-field frame, see \cite{SusanneEva}. 

Notice also that the class of Erlang memory kernels is dense in $ L^1 ( \R_+), $ see e.g.\ \cite{kammler}. Therefore, any Hawkes process $N$ having general integrable memory kernel $h$ can be approximated by a sequence of Hawkes processes $ N^{(n)} $ having Erlang memory kernel $h^{(n)}  $  such that $ \| h^{(n)} - h \|_{L^1 ( \R_+) } \to 0$ as $ n \to \infty $ and 
\begin{equation}\label{eq:approx}
E \| \! | N - N^{(n)}  \|\! |_t	\le C_{T} \int_{0}^{t} | h^{(n)} - h| (s) ds ,
\end{equation}
for all $t \le T $ (see Lemma (3.4) of \cite{madsevasusanne}), where $ \| \! | N - N^{(n)}  \|\! |_t $ denotes the total variation distance between $ N $ and $N^{(n)} $ on $ [0, t ].$   

Finally, the Markovian representation of Hawkes processes having memory kernels as in Assumption \ref{ass:h} in terms of the PDMP \eqref{def:PDMP}-\eqref{eq:gen} has two advantages. The first advantage is that 
stability properties and the longtime behavior of such Hawkes processes can be studied via the well-established theory of PDMPs. 
Since it is straightforward to simulate the PDMP \eqref{def:PDMP}-\eqref{eq:gen} (see Section \ref{sec:Simulations}), one can also simulate Hawkes processes with memory kernels given by sum of Erlang kernels by using this representation.
This is the second advantage.

In the next section we discuss stability properties of the associated PDMP \eqref{def:PDMP}-\eqref{eq:gen} with random jump heights. A simulation algorithm for this PDMP will be presented in Section \ref{sec:Simulations}.

\section{Long-time behavior of the associated Markovian cascade with random jump heights}
\label{sec:long_time_term}

In this section we consider the Markov process $X=(X_t)_{t\geq 0}$ taking values in $\R^{\kappa}$ with (possibly) random jump heights. Its generator is given for any smooth and bounded function $g:\R^{\kappa }\mapsto \R$  by
\begin{equation}
\label{def:geranator}
\mathcal{L}g(x)=\langle F(x),\nabla g(x)\rangle +f\big(\sum_{i=1}^L x^{(i, 0)}\big)\int\big(g(x+\sum_{ i=1}^L c_i e_{(i,n_i) })-g(x)\big)G(dc_1,  \ldots , dc_L),
\end{equation}
where $F:\R^{\kappa }\mapsto \R^{\kappa }$ is the vector field associated to the system \eqref{def:sysODEs} and $G(dc_1, \ldots , dc_L )$ is a probability measure on $\R^L .$ 

\begin{assumption}\label{ass:3}
The probability measure $G$ on $\R^L$ has finite first moments, i.e., 
\begin{equation}\label{eq:finitefirstmoment}
 \int \sum_{ i=1}^L |c_i| G(dc_1, \ldots , d c_L ) < \infty .
\end{equation}
\end{assumption}

The above process is well-defined under Assumptions \ref{ass:Lip} and \ref{ass:3}, as shows the following proposition.

\begin{proposition}\label{prop:noexplosion}
Assume Assumptions \ref{ass:Lip} and \ref{ass:3}. Let $N=(N_t)_{t \geq 0}$ be the counting process associated to the jumps of the Markov Process $X= (X_t)_{t \geq 0}$ having generator given by \eqref{def:geranator}, starting from $x_0\in\R^{\kappa}$. Then $N$ has $P_{x_0}$-almost surely a finite number of jumps on each interval $ [s, t ], 0 \le s < t < \infty .$ 
\end{proposition}
The proof of this proposition is analogous to the proof of Proposition \ref{thm:1}.

\subsection{A Foster-Lyapunov condition}
We start showing that there exists a compact set $K$ of $\R^{\kappa}$ such that the process $X=(X_t)_{t\geq 0}$ possessing the generator defined in \eqref{def:geranator} visits $K$ infinitely often almost surely. Let $n= \max_{1\leq i\leq L} n_i $ and $\alpha= \min_{1\leq i\leq L} \alpha_i.$ In what follows, we write $ 0^{\kappa } $ to denote the vector in $\R^\kappa $ having all coordinates equal to $0.$ 

\begin{proposition}
\label{Prop:driftconditio}
Suppose Assumptions \ref{ass:Lip} and \ref{ass:3}. Let  $\mathcal{L}$ be the generator defined in \eqref{def:geranator} and consider the function $V:\R^{\kappa }\mapsto \R_+$ defined by
\begin{equation}
V(x)=1+\sum_{i=1}^L \sum_{k=0}^{n_i}\frac{b(k+1)}{\alpha_i^k}|x^{(i,k)}|,
\end{equation}
where $b:\{0, 1,\ldots n+1\}\to \R_{+}$ is a strictly  increasing function. 
If $f$ is not bounded but only Lipschitz continuous, we suppose moreover that 
\begin{equation}
\label{ass.1}
\|f\|_{Lip}  \left(  \int \sum_{i=1}^L   \frac{1}{\alpha^{n_i}}|c_i| G (dc_1, \ldots, dc_L) \right) < \alpha 
\end{equation}
and choose the function $b$ so that
\begin{equation}\label{eq:b}
 \frac{b(n+1)}{b(1)} \|f\|_{Lip} \left(  \int \sum_{i=1}^L   \frac{1}{\alpha^{n_i}}|c_i| G (dc_1,\ldots, dc_L) \right)  < \alpha  .
\end{equation}

Then there exist positive constants $\lambda$, $\beta$ and $R$ such that the following Foster-Lyapunov type drift condition holds
\begin{equation}
\label{driftcond}
\mathcal{L}V(x)\leq -\lambda V(x)+\beta 1_{K}(x),
\end{equation}
where $K=\bar{B}_R(0^{\kappa})$ is the (closed) ball of center $0^{\kappa}$ and radius $R$.
\end{proposition}  
\begin{remark}
If $ \alpha_i = \alpha $ for all $1 \le i \le L $ and $sg (c_i ) = sg (c_j ) $ (where $sg(u)$ is the sign of $u\in \R$) for all $ i \neq j ,$ then condition \eqref{ass.1} is equivalent to the sub-criticality condition \eqref{cond:Bremaud} required in Theorem 1 of \cite{bm}. For values of $\alpha = \min_{1\leq i\leq L} \alpha_i \geq 1$, we could have taken the simpler Lyapunov function $V(x)=1+\sum_{i=1}^L \sum_{k=0}^n|x^{(i, k)}|$ .  
\end{remark}

\begin{proof}
Indeed, one immediately verifies that
$$
\mathcal{L}V(x)=A(x)+B(x),
$$
where  
$$
A(x)=\sum_{i=1}^L \left(\sum_{k=0}^{n_i-1}\frac{b(k+1)}{\alpha_i^k}\big(x^{(i, k+1)}-\alpha_i x^{(i, k)}\big)sg\big(x^{(i, k)}\big) -\frac{b(n_i+1)}{\alpha_i^{n_i-1}}|x^{(i, n)}|\right)$$
and
$$
B(x)= f\big(\sum_{i=1}^L x^{(i, 0)}\big) \int \sum_{i=1}^L \frac{b(n_i+1)}{\alpha_i^{n_i}}\big(|x^{(i,n_i)}+c_i|-|x^{(i, n_i)}|\big)G(dc_1,\ldots, dc_L) .
$$

Defining $b_*=\min\{b(k+1)-b(k),0\leq k\leq n\}$ and $r=\alpha b_{*}(b(n+1))^{-1}$,
it is also straightforward to check that
\begin{multline}
\label{prop1:ineq1}
A(x)\leq - b(1) \sum_{i=1}^L  \alpha_i |x^{(i, 0)}| -\alpha \sum_{i=1}^L \sum_{k=1}^{n_i }\frac{1}{\alpha_i^{k}} \big(b(k+1)-b(k)\big)|x^{(i, k)}| \\
\leq- \alpha  \sum_{i=1}^L  b(1)  |x^{(i, 0)}| - r \sum_{i=1}^L \sum_{k=1}^{n_i } \frac{ b(k+1)}{\alpha_i^k} |x^{(i, k )}| .
\end{multline}

Suppose first that is $f$ is bounded by $f^* .$ In this case, one can easily verify that 
$$
B(x)\leq f^{*}b(n+1)\int \sum_{i=1}^L \alpha_i^{-n_i}|c_i|G(dc_1,\ldots, dc_L).
$$ Since $ r < \alpha , $ it follows from the above estimates that 
$$
\mathcal{L}V(x)\leq  -rV(x)+p,
$$
where $p=r+f^{*}b(n+1)\int \sum_{i=1}^L \alpha_i^{-n_i}|c_i|G(dc_1,\ldots, dc_L)$. 

Let $q=1+b(1)(1\wedge \alpha^{-n})R$ and observe that $V(x)\geq q$, for $x\notin K=\bar{B}_R(0^{\kappa})$. Thus, taking any $R$ sufficiently large such that $ p/q < r , $ we deduce that
$$
\mathcal{L}V(x)\leq -\Big(r-\frac{p}{q}\Big)V(x)+p1_{K}(x),
$$
which proves \eqref{driftcond} for bounded jump rates $f$ with $\lambda=r-\frac{p}{q}>0$ and $\beta=p>0$.

Assuming now that $f$ is unbounded and Lipschitz continuous, we have that
$$
B(x)\leq \Big(\|f\|_{Lip} \sum_{i=1}^L |x^{(i, 0)}|+f(0)\Big)b(n+1) \int \sum_{i=1}^L \alpha_i^{-n_i}|c_i|G(dc_1,\ldots, dc_L) ,
$$
which, together with the first inequality in \eqref{prop1:ineq1}, implies that
$$
\mathcal{L}V(x)\leq - d V(x)+d+f(0)b(n+1) \int \sum_{i=1}^L \alpha_i^{-n_i}|c_i|G(dc_1,\ldots, dc_L),
$$  
where $d=\Big(\alpha -\|f\|_{Lip}\frac{b(n+1)}{b(1)}\int  \sum_i \alpha_i^{-n_i} |c_i|G(dc_1,\ldots, dc_L)\Big)\wedge r $ is positive thanks to \eqref{ass.1}. Using the inequality above and proceeding as before we establish also the drift condition \eqref{driftcond} for Lipschitz jump rates.  
\end{proof}

As a corollary of Proposition \ref{Prop:driftconditio}, we obtain exponential moments for the return times to the compact set $K$ appearing in \eqref{driftcond}.

\begin{corollary}\label{cor:1}
Let $K=\bar{B}_R(0^{\kappa})$ and $V$ be as in Proposition \ref{Prop:driftconditio}. Write $T_K=\inf\{t>0:X_t\in K\}.$ Then for all $\eta\leq \lambda$ and $x_0\in\R^{\kappa}$, 
\begin{equation}
E_{x_0}[e^{\eta T_K}]\leq V(x_0).
\end{equation}
\end{corollary}   
The proof of this corollary is classical, see for instance Theorem 6.1 of \cite{dmt}.

\subsection{Wasserstein contraction for Lipschitz jump rates}
Throughout this section we suppose that the jump rate $f$ is Lipschitz continuous. In this case, we are able to prove the exponential convergence to equilibrium in Wasserstein distance,  under the sub-criticality condition \eqref{ass.1}. 

More precisely, in the sequel, for any $ x \in \R^{\kappa }, $ we will write $\|x\|_1 = \sum_{i=1}^L \sum_{k=0}^{n_i} |x^{(i,k)}| .$ Let $\mu $ and $ \nu $ be two probability measures on $\R^{\kappa}.$ We call coupling of $ \mu $ and $ \nu $ any probability measure on $ \R^{ \kappa} \times \R^{\kappa}$ whose marginals are $ \mu $ and $ \nu, $ and we denote by $ \Gamma ( \mu, \nu ) $ the set of all such couplings. The Wasserstein distance between $\mu $ and $\nu$ is defined by 
\begin{equation}\label{eq:Wasserstein}
W_1 ( \mu , \nu ) = \inf \left\{ \int_{ \R^{\kappa }} \int_{\R^{\kappa }} \| x - y\|_1 \gamma (dx, dy ) , \gamma \in   \Gamma ( \mu, \nu ) \right\} .
\end{equation}

In the following, we write $(P_t)_{t \geq 0}$ for the transition semigroup of the process $X$ with generator \eqref{def:geranator}. Recall that $ A = \max_{1\leq i\leq L} \alpha_i $, $\alpha = \min_{1\leq i\leq L} \alpha_i$ and  
$n= \max_{1\leq i\leq L} n_i $. The following theorem states the exponential rate of convergence to equilibrium of the process with respect to the Wasserstein distance.

\begin{thm}\label{theo:wassersteincontraction}
Suppose $f$ is Lipschitz continuous, assume condition \eqref{ass.1} and choose the function $b$ as in \eqref{eq:b}. 

1. Then, for any choice of probability measures $\mu $ and $\nu $ on $ {\cal B}(\R^{\kappa}), $ 
\begin{equation}\label{eq:W1control}
 W_1 ( \mu P_t, \nu P_t) \le C  e^{- d t } W_1 ( \mu , \nu ) , 
\end{equation} 
where 
$$ C = \frac{ A^n \vee 1}{1 \wedge \alpha^{n+1}} \frac{b(n+1)}{b(1) } $$ 
and $d=\Big(\alpha -\|f\|_{Lip} \frac{b(n+1)}{b(1)} \int \sum_{i=1}^L \alpha_i^{-n_i}|c_i|G(dc_1,\ldots, dc_L)\Big)\wedge \Big( \alpha b_{*}(b(n+1))^{-1}\Big), $ with $b_*=\min\{b(k+1)-b(k),0\leq k\leq n\}.$

2. In particular, there exists a unique invariant probability measure $\pi $ of the process $X$ such that for any probability measure  $\nu $ on $ {\cal B}(\R^{\kappa}), $ 
$$  W_1 ( \pi , \nu P_t) \le C  e^{- d t } W_1 ( \pi  , \nu ) . $$
\end{thm}

\begin{proof}
The assertion of point 1. follows from a standard Wasserstein coupling. More precisely, denote by $ (X, \tilde X) $ the Markov processes taking values in $ \R^{\kappa} \times \R^{\kappa} $ having the infinitesimal generator defined for any
smooth test function $ \varphi ( x, y ) : \R^{\kappa} \times \R^{ \kappa } \to \R  $ by  
\begin{multline}\label{eq:generatorcoupling}
\mathcal{L}_2  \varphi ( x, y) = \langle F(x), \nabla_x\varphi(x,y) \rangle + \langle F(y), \nabla_y\varphi(x,y) \rangle
\\
+  f\big( \sum_i x^{(i, 0) }\big) \wedge f\big( \sum_i y^{(i,0) }\big) \int G(dc_1,\ldots, dc_L)  \left[ \varphi ( x +\sum_i c_i e_{(i,n_i) }, y +\sum_i c_i e_{(i,n_i)} ) - \varphi ( x, y ) \right]  \\
+  \left( f \big( \sum_i x^{(i,0)}\big) -  f \big(\sum_i y^{(i,0)}\big)\right)_+ \int G(dc_1,\ldots, dc_L) \left[ \varphi ( x +\sum_i c_i e_{(i,n_i)} , y ) - \varphi ( x, y ) \right] \\
+   \left( f \big(\sum_i  y^{(i,0)}\big) -  f \big(\sum_i  x^{(i,0)}\big)\right)_+ \int G(dc_1,\ldots, dc_L)  \left[ \varphi ( x, y + \sum_i c_ie_{(i,n_i)} ) - \varphi ( x, y ) \right],
\end{multline}
where $F:\R^{\kappa}\mapsto \R^{\kappa}$ is the vector field associated to the system \eqref{def:sysODEs}.

This is the usual coupling which consists of making the two processes jump together as much as possible. Define 
$$ H ( x, y )= \sum_{i=1}^L \sum_{k=0}^{n_i} \frac{ b(k+1)}{\alpha_i^k }   |x^{(i,k)}- y^{(i,k)  }| .$$
Then an analogous calculus as the one used in the proof of Proposition \ref{Prop:driftconditio} yields
$$ \mathcal{L}_2 H (x, y ) \le - d H(x, y ) $$
implying that 
$$ E_{x, y } H( X_t, \tilde X_t) \le H(x, y ) e^{- d t}.$$ 
Observing that 
$$ \| x - y \|_1 \le \frac{ A^n \vee 1 }{b(1) } H(x, y ) ; \; H(x, y ) \le \frac{b(n+1)}{1 \wedge \alpha^{n+1}} \| x- y \|_1 $$
we conclude the proof of item 1. 

To prove item 2., let $ \mu_n := \mu P_n$ for any probability measure $\mu $ on $ {\cal B}(\R^{\kappa}).$ Observe that $ W_1 ( \mu_{n+m}, \mu_m ) \le C e^{- dm } W_1 ( \mu_n , \mu ) ,$ implying that $ (\mu_n)_n $ is Cauchy and thus, by the completeness of the space of all probability measures on $ (\R^\kappa , {\cal B} ( \R^\kappa ) ) ,$ endowed with the metric induced by $W_1 $ (see e.g.  \cite{Rachev} or \cite{bolley}), convergent to some limit measure $ \mu_\infty .$ This limit measure must be invariant. Indeed we have  $ W_1 ( \mu_m P_t, \mu_m ) \le C e^{- d m } W_1 ( \mu P_t, \mu ) .$ But 
\begin{eqnarray*}
 |  W_1 ( \mu_m P_t, \mu_m )  -  W_1 ( \mu_\infty  P_t, \mu_\infty ) | &\le& W_1 ( \mu_m P_t, \mu_\infty P_t  ) + W_1 ( \mu_\infty  , \mu_m )\\
 &\le &[C e^{-d t} +1 ]  W_1 ( \mu_\infty  , \mu_m ) \to 0 
\end{eqnarray*} 
as $m \to \infty , $ where we have used \eqref{eq:W1control} once again to obtain the second inequality. As a consequence, $W_1 ( \mu_\infty P_t , \mu_\infty ) = 0, $ implying that $\mu_\infty $ is the (necessarily unique) invariant measure. This concludes our proof. 
\end{proof}

\begin{remark}\label{rem:3}
Under the conditions of Theorem \ref{theo:wassersteincontraction}, write $ \bar N_t $ for the stationary version of the non-linear Hawkes process having intensity \eqref{eq:intensity}; that is, $\bar N_t $ has intensity $ \bar \lambda_t  =  f \big(\sum_{i=1}^L \bar X^{(i, 0)}_ t \big), $ where $ \bar X $ is the stationary process evolving according to \eqref{eq:gen}. Let moreover $ N_t $ be the non-linear Hawkes process with intensity $ \lambda_t = f \big(\sum_{i=1}^L  X^{(i, 0)}_ t \big) $ starting from some fixed initial condition $X_0 = x_0 \in \R^\kappa .$ Then the Lipschitz-continuity of $f$ together with \eqref{eq:W1control} imply that $ E  \int_0^\infty | \lambda_t - \bar \lambda_t | dt < \infty .$ It is then straightforward to deduce from this by standard coupling arguments, as explained e.g.\  the proof of Theorem 1 in \cite{bm}, that $ N $ and $ \bar N$ couple almost surely  in finite time; that is, there exists $ T > 0 $ such that for all $ t \geq 0, $ $ N_{T+t}- N_T = \bar N_{T+t} - \bar N_t,$ meaning that $ N$ and $ \bar N$ have the same jump times after time $T.$ This is what is called {\it stability in variation} in \cite{bm} (see their Definition 1). Therefore, our Theorem \ref{theo:wassersteincontraction} implies Theorem 1 of \cite{bm}. 
\end{remark}

In the next section we prove a stronger result, showing that the process is even recurrent in the sense of Harris.

\subsection{Harris recurrence}
In this section, we use the regeneration method based on Nummelin splitting to show that $X$ is recurrent in the sense of Harris having a unique invariant probability measure $\pi .$ We recall (e.g. from \cite{adr}) that 
\begin{definition}
\label{def:harris_recurrent}
The process $(X_t)_{t\geq 0}$ is said to be \textit{recurrent in the sense of Harris} if there exists a sigma-finite measure $m $ on ${\cal B}(\R^{\kappa})$  such that $m (A)>0$  implies that  for all $\ x\in \R^{\kappa}$, $P_x-$almost surely, 
$$
\limsup_{t\to \infty} 1_A(X_t)=1 .
$$
\end{definition}
By \cite{adr},  Harris recurrence of $X$ implies in particular the existence of a unique invariant measure (which is sigma-finite but does not need to be finite) $\pi$ such that  the above property holds with $\pi $ in place of $m$. $X$ is called positive Harris recurrent if $ \pi ( \R^{\kappa } ) < \infty .$ We have the following

\begin{thm}
\label{harrisrecofY_X}
Suppose that $f$ is bounded or Lipschitz continuous satisfying \eqref{ass.1}.  Suppose moreover that Assumption \ref{ass:3} holds and that $ G( dc_1, \ldots , d c_L) = \prod_{i=1}^L G_i (d c_i ) $ for probability measures $ G_i $ on $ ( \R, {\cal B} (\R) ) $ satisfying supp $(G_i) \cap \{0\}^c \neq \emptyset , $ for all $ 1 \le i \le L.$ Finally, suppose that $f$ is lower bounded. 

1. Then $ (X_t)_{t \geq 0} $ is positive Harris recurrent with unique invariant measure $ \pi (dx) . $ \\

2. Let $ \bar X_t $ be a stationary version of the process  and suppose that $ (X_t)_{t \geq 0} $ starts from $ X_0 = x_0 \in \R^\kappa ,$ both evolving according to \eqref{eq:gen}.  Then $\bar X$ and $ X$ couple almost surely in finite time; that is, it is possible to construct them on the same probability space such that there exists $ \tau_c < \infty $ almost surely satisfying 
\begin{equation}\label{eq:couplingofx}
t \geq \tau_c \mbox{  implies that   } X_t = \bar X_t \; \mbox{ and } \;  P ( \tau_c > t ) \le C (p )  V(x) t^{-p} , 
\end{equation}
for every $ p \geq 1, $ where $C(p) $ is a constant depending on $p.$ 
\end{thm}

\begin{remark}
In particular, using the notation of Remark \ref{rem:3} above, \eqref{eq:couplingofx}
implies that $ \lambda_t = \bar \lambda_t $ for all $t \geq \tau_c $ meaning that $\bar N$ and $N$ couple as well. As a consequence, our Theorem \ref{harrisrecofY_X} is a refinement of the results of Theorem 1 and Theorem 2 in \cite{bm} -- however at the prize of imposing a lower bound on $f.$ 
\end{remark}

The proof of this theorem uses the regeneration technique based on Nummelin splitting. It is well known that it is easier to implement this method in the frame of discrete time Markov chains rather then Markov processes in continuous  time -- although some effort has been spent to introduce regeneration times in a continuous time framework, see e.g.\ \cite{dashaeva}. Therefore, we start by showing that the sampled chain $ (Y_k)_{k \geq 0} = (X_{k T})_{ k \geq 0 } , $ for some fixed $T > 0 ,$ is positive Harris recurrent. 

We recall that the chain $(Y_k)_{k\geq 0}$ is said to be \textit{recurrent in the sense of Harris} with invariant measure $\pi $ on ${\cal B}(\R^{\kappa})$ if whenever $\pi (A)>0$, we have, for all $ x\in \R^{\kappa} ,$  ${P}_x-$almost surely, 
$
  \limsup_{k\to \infty} 1_A(Y_k)=1.
$
Obviously, Harris recurrence of the chain $ (Y_k)_{k \geq 0} $ implies the Harris recurrence of the process $ X, $ and the invariant probability measures of both processes coincide (if they exist).  

The rest of this section is devoted to prove that the sampled chain $(Y_k)_{k \geq 0} $ is Harris which follows from the following Doeblin type lower bound. Recall that $(P_t)_{t \geq 0}$ denotes the transition semigroup of the process $X,$ therefore, $P_T $ is the transition operator of the sampled chain $  (Y_k)_{k \geq 0} .$ 

\begin{thm}\label{thm:Doeblin}
Assume the assumptions of Theorem \ref{harrisrecofY_X}. For all $x^*  \in \R^{\kappa}, $ there exist $R > 0 , $ an open set $ I \subset \R^\kappa $ and a constant $\beta \in (0, 1), $ depending on $I,R, L, n_1, \ldots, n_L, \alpha_1,\ldots, \alpha_L$ and $f$, such that 
\begin{equation}\label{doblinminorization}
P_{LT} ( x, dy ) \geq \beta 1_C  (x) \nu ( dy) ,
\end{equation} 
where $ C = B_R ( x^* ) $ is the (open) ball of radius $R$ centered at $x^* ,$ and where $ \nu $ is the uniform probability measure on $  I.$ 
\end{thm}   

\begin{proof}
{\bf Part I. $L=1.$}\\
We start by proving the result in the case $L=1, c_1=c,\alpha_1=\alpha$ and $n_1=n$, that is, $h( t) = c e^{ - \alpha t} \frac{t^n }{n!}.$ The corresponding Markov process is then given by $ X_t = (X_t^{(0)}, \ldots, X_t^{(n)}) $ taking values in $\R^{n+1}.$ Clearly, for all $ A \in {\cal B} ( \R^{n+1}),$ 
$$ P_T ( x, A) \geq E_x ( 1_A ( X_T), N_T = n+1) .$$ 

Recall the definition of the flow in \eqref{eq:flow}. On the event $ \{ N_T = n+1 \}, $ starting from $X_0 = x , $ we first let the flow evolve starting from $ x  $ up to some first jump time $t_1 .$ At that jump time we choose an associated jump height $c_1.$ We then successively choose the following inter-jump waiting times $t_2, \ldots, t_{n+1}$ under the constraint $ t_1 + \ldots + t_{n+1} < T $ and the associated jump heights $c_2, \ldots , c_{n+1}.$ Write $ s_1= T- t_1, s_2= T - (t_1 + t_2) , \ldots, s_{n+1}  = T - (t_1 + \ldots + t_{ n+1}) .$ 

Conditionally on $X_0=x$, the successive choices of ${\bf \underline{c}}=(c_1, \ldots, c_{n+1}) $ and ${\bf\underline{s}}=(s_1, \ldots , s_{n+1})$ as above, the position of $ X_T $ is given by 
\begin{equation}\label{def:flowofZ}
\gamma (x,{\bf \underline{c}} , {\bf \underline{s}})=\varphi_T  ( x  )+c_1 e^{- \alpha s_1}v({s_1})+\ldots + c_{n+1} e^{- \alpha s_{n+1}}v({s_{n+1}}), 
\end{equation}
where for each $1\leq k\leq n+1$,
\begin{equation}
v({s_k})=
\left( \begin{array}{c}
\frac{s_{k}^n}{n!} \\
\frac{s_{k}^{n-1}}{(n-1)!} \\
\vdots \\
s_{k}\\
1 
\end{array}
\right).
\end{equation}
We omitted the dependence on $T$ of the map $\gamma ( x,{\bf \underline{c}} , {\bf \underline{s}})$ since we keep the value $T>0$ fixed once for all and work with sequences ${\bf \underline{s}}$ satisfying the constraints $0< s_{n+1}< s_{n}< \ldots < s_1 < T$. 
Finally, in what follows we shall write, for any fixed pair $ (x ,{\bf \underline{c}}),$  
$$ \gamma_{ (x ,{\bf \underline{c} } ) } :    {\bf \underline{s}} \mapsto \gamma (x,{\bf \underline{c}} , {\bf \underline{s}}). $$

We will use the jump noise which is created by the $n+1$ jumps, i.e., we will use a change of variables on the account of $s_1, \ldots, s_{n+1} .$ 
Therefore, in what follows we write 
$$
\frac{\partial \gamma_{x, {\bf \underline{c}}}({\bf \underline{s}})}{\partial {\bf \underline{s}}}=\Big[\frac{\partial \gamma_{x, {\bf \underline{c}}}({\bf \underline{s}})}{\partial s_1},\ldots, \frac{\partial \gamma_{x, {\bf \underline{c}}}({\bf \underline{s}})}{\partial s_{n+1}} \Big]
$$ 
to denote the Jacobian matrix of the the map $\underline{\bf s}\mapsto \gamma_{x, {\bf \underline{c}}}({\bf \underline{s}})$.  
This matrix does not depend on the initial position $x$ nor on the first jump height $c_0$. Indeed, one easily finds that
$$
\frac{\partial \gamma_{x, {\bf \underline{c}}}({\bf \underline{s}})}{\partial {\bf \underline{s}}}=\Big[C^{(1)},\ldots ,C^{(n+1)} \Big],
$$ 
where for each $1\leq k\leq n+1$, $C^{(k)}$ is a column vector given by 
$$
C^{(k)}=c_ke^{-\alpha s_k}
\left( \begin{array}{c}
\frac{s_{k}^{n-1}}{(n-1)!}-\alpha \frac{s_{k}^n}{n!}\\
\frac{s_{k}^{n-2}}{(n-2)!}-\alpha \frac{s_{k}^{n-1}}{(n-1)!} \\
\vdots \\
1-\alpha s_{k}\\
-\alpha 
\end{array}
\right).
$$
As a consequence the determinant of $\frac{\partial \gamma_{x, {\bf \underline{c}}}({\bf \underline{s}})}{\partial {\bf \underline{s}}}$ is given by
\begin{equation}\label{eq:determinant}
\mbox{det}\Big(C^{(1)},\ldots, C^{(n+1)}\Big)=(-1)^{n+1}\alpha \prod_{k=1}^{n+1}c_ke^{-\alpha s_k}\mbox{det}\Big(\tilde{C}^{(1)},\ldots, \tilde{C}^{(n+1)}\Big),
\end{equation}
where for each $1\leq k\leq n+1$, $\tilde{C}^{(k)}$ is a column vector given by 
$$
\tilde{C}^{(k)}=
\left( \begin{array}{c}
\alpha \frac{s_{k}^n}{n!}-\frac{s_{k}^{n-1}}{(n-1)!} \\
\alpha \frac{s_{k}^{n-1}}{(n-1)!}-\frac{s_{k}^{n-2}}{(n-2)!} \\
\vdots \\
\alpha s_{k}-1\\
1 
\end{array}
\right).
$$
Thus the invertibility of the matrix $\frac{\partial \gamma_{x, {\bf \underline{c}}}({\bf \underline{s}})}{\partial {\bf \underline{s}}}$ follows from the invertibility of the matrix $J=[\tilde{C}^{(1)},\ldots , \tilde{C}^{(n+1)}]$. In the sequel, for each $1\leq k\leq n+1$, let $r_k$ denote the $k$-th row of $J$. 
By replacing successively (bottom-up) $r_i$ by $\alpha^{-1}(r_{i}+r_{i-1})(n+1-i)!$, we deduce that $J$ is equivalent to the Vandermonde matrix
$$
\begin{pmatrix}
s_1^{n} & s_2^{n} & \ldots & s_{n+1}^{n} \\
s_1^{n-1} & s_2^{n-1} & \ldots & s_{n+1}^{n-1} \\
\vdots \\
s_1 & s_2 & \ldots & s_{n+1} \\
1 & 1 & \ldots & 1 
\end{pmatrix},
$$
which is know to be invertible if and only if $0<s_{n+1}< s_n <\ldots <s_{1}$. In conclusion, we have just shown that for any $x\in \R^{n+1}$, any choice of ${\bf \underline{c}}$ having non null coordinates,  
the Jacobian of the map $\underline{s}\mapsto \gamma_{x,{\bf \underline{c}}}({\bf \underline{s}})$ is invertible at any ${\bf \underline{s}}$ such that $0<s_{n+1}< s_n <\ldots <s_{1}$. 

It will be proved now that this uniform invertibility of the Jacobian matrix of the map $\underline{\bf s}\mapsto \gamma_{x,{\bf \underline{c}}}({\bf \underline{s}})$ implies  inequality \eqref{doblinminorization}. For that sake, we shall also need the following notation. For each triple $(x, {\bf \underline{c}}, {\bf \underline{s}})$, we write $x_0 = x , $ $x_1=\varphi_{T-s_{1}}(x )+c_1e_{n+1}$ \footnote{$e_{n+1}$ denotes the $n+1-$st unit vector in $ \R^{n+1} $} and then recursively $x_k=\varphi_{s_{k-1}-s_{k}}(x_{k-1})+c_ke_{n+1}$ for $2\leq k\leq n+1$. The sequence $x_1,\ldots x_{n+1}$ corresponds to the positions right after successive jumps, starting from the initial location $x\in \R^{n+1},$ induced by the heights ${\bf \underline{c}}$ and the inter-jump waiting times $T-s_1,s_1-s_2, \ldots s_n-s_{n+1} $ which are determined by ${\bf \underline{s}}$. 

Introduce now for each $x\in \R^{n+1}$ and $t\geq 0,$ 
\begin{equation}
\label{def:suvival_rate}
e(x,t)= \exp\Big\{-\int_{0}^tf\big(\varphi_s^{(0)}(x)\big)ds\Big\}
\end{equation}
and define  for each triple  
$(x, {\bf \underline{c}}, {\bf \underline{s}})$ (here we set $s_0=T$), 
\begin{equation}\label{densityat_x_c_s}
q_{x,{\bf \underline{c}}}({\bf \underline{s}})=\left( \prod_{k=0}^{n}f(\varphi_{s_k-s_{k+1}}^{(0)}(x_k))e(x_k,s_k-s_{k+1}) \right) e(x_{n+1},s_{n+1}).
\end{equation}

Since $f$ is bounded away from 0 and from the definition of $e(\cdot, \cdot),$ we deduce that for any triple $(x^*,{\bf \underline{c}}^*,{\bf \underline{s}^*})$ there are neighborhoods $W_{\bf \underline{s}^*}$, $U_{x^*}$ and $V_{\bf \underline{c}^*}$ of ${\bf \underline{s}}^*, $ ${x^*}$ and ${\bf \underline{c}^*}$ respectively such that  
\begin{equation}
\label{lowerboundq}
\inf_{(x,{\bf \underline{c}}, {\bf \underline{s}})\in U_{x^*}\times V_{\bf \underline{c}^*} \times W_{\bf \underline{s}}^*} q_{x,{\bf \underline{c}}}({\bf \underline{s}})>0. 
\end{equation}   
Let us now fix a triple $(x^*,{\bf \underline{c}}^{*},{\bf \underline{s}}^{*})$ such that the matrix $\frac{\partial \gamma_{x^*, {\bf \underline{c}}^*}({\bf \underline{s}}^*)}{\partial {\bf \underline{s}}}$ is invertible.
Recall that by \eqref{eq:determinant}, the vector ${\bf \underline{c}}^{*}$ must have all coordinates non-null. 
By Lemma 6.2 of \cite{benaim2015}, there exist an open neighborhood $J_{x^*,{\bf \underline{c}}^*}= {B}_{R}(x^*)\times  {B}_{R} ({\bf \underline{c}}^*) $ of the pair $(x^*,{\bf \underline{c}}^*),$ an open set $I\subset \R^{n+1}$, and for any pair $(x,{\bf \underline{c}})\in J_{x^*,{\bf \underline{c}}^*},$ an open set $W_{x,{\bf \underline{c}}}$ such that
$$
\tilde{\gamma}_{x,{\bf \underline{c}}}({\bf \underline{s}}):\left\lbrace
\begin{array}{c}
W_{x,{\bf \underline{c}}} \to I \\
{\bf \underline{s}} \mapsto \gamma_{x,{\bf \underline{c}}}({\bf \underline{s}}),
\end{array}
\right.
$$
is a diffeomorphism, with $W_{x,{\bf \underline{c}}}\subset W_{\bf \underline{s}}^*$ 
and also
\begin{equation}
\label{lowerbounddet}
\inf_{x,{\bf \underline{c}}\in J_{x^*,{\bf \underline{c}}^*}} \inf_{\underline{s}\in W_{x,{\bf \underline{c}}}} \Big|\mbox{det}\Big(\frac{\partial \gamma_{x, {\bf \underline{c}}}({\bf \underline{s}})}{\partial {\bf \underline{s}}}\Big)^{-1}\Big|>0.
\end{equation}
Reducing (if necessary)   
$R$, we may assume also that $J_{x^*,{\bf \underline{c}}^*}\subset U_{x^*}\times V_{\bf \underline{c}^*}$. Thus we have that by  \eqref{lowerboundq} and \eqref{lowerbounddet}, 
\begin{equation}
\label{lowerbounddet2}
\inf_{x,{\bf \underline{c}}\in J_{x^*,{\bf \underline{c}}^*}} \inf_{\underline{s} \in W_{x,{\bf \underline{c}}}} q_{x,{\bf \underline{c}}}({\bf \underline{s}}) \Big|\mbox{det}\Big(\frac{\partial \gamma_{x, {\bf \underline{c}}}({\bf \underline{s}})}{\partial {\bf \underline{s}}}\Big)^{-1}\Big|>0.
\end{equation}
Since $supp(G)\cap\{0\}^c\neq\emptyset$ 
there exists an interval $(a,b)$ such that $0\notin (a,b)$ and $G((a,b))>0$.
Thus, by taking ${\bf \underline{c}}^{*}=((a+b)/2, \ldots, (a+b)/2)$,  we have (reducing $R$ again if necessary) that for $1\leq k\leq n$, 
$G(({\underline{c}^*}^{(k)}-R,{\underline{c}^*}^{(k)}+R))>0$ which together with \eqref{lowerbounddet2} implies
\begin{equation}\label{eq:tildebeta}
\tilde \beta=\Big(\prod_{k=1}^nG(({\underline{c}^*}^{(k)}-R,{\underline{c}^*}^{(k)}+R))\Big)\inf_{x,{\bf \underline{c}}\in J_{x^*,{\bf \underline{c}}^*}} \inf_{\underline{s}\in W_{x,{\bf \underline{c}}}} q_{x,{\bf \underline{c}}}({\bf \underline{s}}) \Big|\mbox{det}\Big(\frac{\partial \gamma_{x, {\bf \underline{c}}}({\bf \underline{s}})}{\partial {\bf \underline{s}}}\Big)^{-1}\Big|>0.
\end{equation}
Finally, we have for any measurable
$A \in \mathbb{B}(\R^{n+1})$ and $x\in  {B}_R (x^*), $ using the change of variables $y=\tilde{\gamma}_{x,{\bf \underline{c}}}({\bf \underline{s}})$,
\begin{eqnarray}\label{eq:tobe}
E_x ( 1_A ( X_T, N_T = n+1) &\geq &\int_{ {B}_{R}({\bf \underline{c}}^*)}\underline{G}(d{\bf \underline{c}}) \int_{W_{x,{\bf \underline{c}}}} q_{x,{\bf \underline{c}}}({\bf \underline{s}})1_{A}(\gamma_{x, {\bf \underline{c}}}({\bf \underline{s}}))ds_{1}\ldots ds_{n+1} \nonumber\\
&\geq & \tilde \beta \int_{I\cap A}dy_1\ldots dy_{n+1}=\beta \nu(A),
\end{eqnarray}
where $\underline{G}(d{\bf \underline{c}})=G(dc_1)\ldots G(dc_{n+1})$ and $ \beta = \tilde \beta \nu ( I) , $ establishing the desired result in case $L=1.$

The proof of the general case $L > 1 $ follows the same strategy and is given in the Appendix. 
\end{proof}

We are now able to conclude the proof of Theorem \ref{harrisrecofY_X}.

\begin{proof}[Proof of Theorem \ref{harrisrecofY_X}]
1) By Corollary \ref{cor:1}, we know that $X$ comes back to $K $ infinitely often almost surely. Moreover, $ \sup_{x \in K} | \varphi_t (x) | \to 0 $ as $t \to \infty ,$ by the explicit form of the flow in \eqref{eq:flow}. Therefore, for any $ \varepsilon > 0  $ there exists $t^* $ such that $\varphi_t (x) \in B_{\varepsilon } ( 0) $ for all $t \geq t^* , $ for all $ x \in K.$ Since $f$ is bounded on $ \tilde K := \{ \varphi_t (x) : t \geq 0, x \in K \} ,$ we have $ \inf_{x \in K} P_x ( T_1 > t_* + 2 T ) > 0 .$  This implies that 
$$ \inf_{x \in K} P ( \mbox{ the trajectory of } X \in B_{\varepsilon } ( 0) \mbox{ during a time interval of length $> T$ } | X_0 =x ) > 0 , $$ 
and therefore, using a conditional version of the Borel-Cantelli lemma, $(Y_k)_{k \in \N}  $ visits $ B_{\varepsilon } ( 0)$ infinitely often almost surely. \\
2) Applying the result of Theorem \ref{thm:Doeblin} with $x^* = 0 $ and $ \varepsilon = R $ and using the standard regeneration technique allows to conclude that $(Y_{Lk})_k $ and therefore $(X_t)_t $ are Harris recurrent. This implies item 1. of the theorem. \\
3) To prove item 2., it is straightforward to show that Proposition \ref{Prop:driftconditio} implies the existence a coupling of $X_t $ and $ \tilde X_t,$ both evolving according to \eqref{eq:gen}, such that for $ T_{K \times K} := \inf \{ t \geq 0 : ( X_t , \tilde  X) \in K \times K \}, $ we have 
$$ E_{x,y} ( e^{\eta  T_{K \times K}})  \le V(x) + V(y).$$
Indeed, if suffices to define the $2\kappa -$dimensional Lyapunov function $ \bar V (x, y ) := V(x) + V(y) $ and to check that \eqref{driftcond} holds for $ \bar {\mathcal L} $ where $ \bar {\mathcal L}$ denotes the generator of the process $ ( X_t, \tilde  X_t ).$  Moreover, \eqref{doblinminorization} can be immediately extended to a lower bound for the joint transition kernel of $ (X_t, \tilde X_t), $ whenever both of them start within the set $ C = B_R (0) .$ Thus $ X $ and $\tilde X$ couple at least with probability $\beta, $ each time they are within $C$ at the same time. The proof that this coupling time has polynomial moments of any order follows then the same arguments as those given in the proof of Proposition 2.15 in \cite{evalast}, implying that 
\begin{equation}\label{eq:last}
 E_{ x, y } [ \tau_c^p ] \le C(p ) [ V(x) + V(y ) ].
\end{equation}  
Finally, Theorem 4.3 of \cite{Meyn93} implies that $ \int V d \pi < \infty $ such that we are able to integrate \eqref{eq:last} against $ \pi ( dx) $ in order to replace $X_t $ by the invariant process $ \bar X_t $ starting from $\bar X_0 \sim \pi .$ This concludes the proof. 
\end{proof}

\section{Simulation Algorithm}
\label{sec:Simulations}
As a consequence of Proposition \ref{thm:1} it follows that 
any Hawkes process possessing memory kernels given by the sum of Erlang kernels can be represented as the counting process associated to the jumps of its  Markovian cascade. Based on this Markovian representation we propose an algorithm (hereafter Algorithm \ref{Alg:unboundedf}) for simulating such Hawkes processes. 

In what follows, for any $x \in \R^{\kappa}, $ we shall write $ \|x\|_\infty = \max \{ |x^{(i,k)}| ,1 \le i \le L,  0 \le k \le n_i \}.$ For a practical implementation of our algorithm the remark below will be important. Recall that $ n= \max_{1\leq i\leq L} n_i$ and $\alpha= \min_{1\leq i\leq L} \alpha_i .$ 
\begin{lemma}
\label{lem:1}
For each $x\in\R^{\kappa},$ let $M(x)=\max\{|\varphi^{(i, 0)}_t(x)| : 1 \le i \le L, t\geq 0\}$ where  $\varphi^{(i, 0)}_t(x)$
is defined in $\eqref{eq:flow}$. Then
\begin{equation}
\label{ineq:Max}
M(x)\leq e\|x\|_\infty  \Big(1\vee \Big(\frac{n}{\alpha e}\Big)^n\Big).
\end{equation}
\end{lemma}

\begin{proof}
Indeed, it follows from \eqref{eq:flow} that for each  $x\in\R^{n+1}$ and $t\geq 0$,
$$
|\varphi^{(0)}_t(x)| \leq \|x\|_\infty e^{-\alpha t}(1+t+\ldots+ t^n/n!)\leq e\|x\|_\infty e^{-\alpha t}(1\vee t^n)\leq e\|x\|_\infty e^{-(\alpha-1)t},
$$
so that if $\alpha>1$, then clearly \eqref{ineq:Max} holds. Now, assume $0<\alpha\leq 1$. Under this assumption, from standard calculus arguments we deduce that $\arg\max \{e^{-\alpha t}(1\vee t^n): t\geq 1\}=n/\alpha.$ This fact and the second inequality above imply the bound in \eqref{ineq:Max} as stated.
\end{proof}

In the sequel, for any rate function $f$ satisfying Assumption \ref{ass:Lip} we define the function $\R^{\kappa }\in x\mapsto f^*(x)$ by 
$$ f^* ( x) =\left\{
\begin{array}{ll}
\max\{f(y):y\in [0, L M(x)]\}, & \mbox{ if } x \in \R_+^{\kappa } \\
\max \{ f(y) : y \in [- L M(x) , 0] \} ,& \mbox{ if } x \in \R_-^{\kappa } \\
\max \{ f(y) : y \in [- L M(x) ,  L M(x) ] \}, & \mbox{ else}
\end{array}
\right\} .$$
Here, $L$ is the number of terms in the sum defining the memory kernel h (recall Assumption \ref{ass:h}).
It follows immediately from Lemma \ref{lem:1} that the function $f^*$ is well-defined, that is  
$f^*(x)$ is finite for all $x\in \R^{\kappa }$.
Let $T_0=0$ and $(T_k)_{k\geq 1}$ denote the sequence of jump times of the Markovian cascade $X$ having generator \eqref{eq:gen}. Observe that the non-explosiveness of $X$ (thanks to Proposition \ref{thm:1}) ensures that the sequence $(T_k)_{k\geq 1}$ is well-defined. Suppose that $X_{T_k}=x$ is given for some $k\in \N$.
Algorithm \ref{Alg:unboundedf} works as follows.
Draw an exponential random variable $\tau$ with parameter $f^*(x)$ and a uniform random variable $U$ on $[0,1]$. If  $U\leq f( \sum_{i=1}^L \varphi^{(i, 0)}_{T_k+\tau}(x))/f^*(x)$, then define the next jump time $T_{k+1}=T_k+\tau$. If not, repeat this procedure starting from $X_{T_k+\tau}=\varphi_{\tau}(x)$ . 
Notice that Algorithm \ref{Alg:unboundedf} is an extension (to our framework) of the classical  thinning algorithm for simulating non-homogeneous Poisson processes. Moreover, it provides an exact simulation of the Markovian cascade $X$ (and consequently of the associated Hawkes process) in the sense that no approximation procedure is required. Its formal definition is given below as a pseudo-code.

\begin{algorithm}[!h]
\caption{Simulation algorithm for the Markovian cascade X}
\label{Alg:unboundedf} 
\begin{algorithmic}[1]
\STATE {\it Input:} bounded or Lipschitz continuous $f$, constants $\alpha_i >0$, $c_i\in\R$ and $T>0$;
and a vector of initial conditions $\big(X^{(i, k )}_0 , 1 \le i \le L, 0 \le k \le n_i \big)=\big(x_0^{(i,k)} , 1 \le i \le L, 0 \le k \le n_i \big)\in \R^{\kappa}.$
\State \textit{Output:} The counting process  $(N_t)_{t\in [0,T]}$.
\State \textit{Initial values:} $x \gets \big(x_0^{(i,k)}, 1 \le i \le L, 0 \le k \le n_i \big) $, $D \gets 0$ and $N_0\gets 0.$ 
\WHILE{$D <T$}
\State $f^*\gets f^*(x)$
\State {draw $\tau\sim \mathcal{E}(f^*)$}
\IF{$\tau\leq T-D  $} 
\State {draw $U\sim U[0,1]$} 
\IF{ $U\leq f( \sum_i \varphi^{(i, 0)}_{D +\tau}(x))/f^*$}
\State $x\gets \varphi_{\tau}(x)+\sum_{i=1}^L c_i e_{(i,n_i)}$
\State $N_t\gets N_D$, for $D \leq t < D+\tau $
\State $N_{ D + \tau} \gets N_D +1 $
\ELSE  
\State $x\gets \varphi_{\tau}(x)$
\State $N_t\gets N_D$, for $D \leq t  \le  D+\tau $ 
\ENDIF
\ELSE \State $N_t\gets N_D$, for $D \leq t \leq T $
\ENDIF 
\State {$D \gets D +\tau$}
\ENDWHILE
\STATE \textbf{Return} $(N_t)_{t\in [0,T]}$.
\end{algorithmic}
\end{algorithm}

Proposition \ref{thm:1} and Lemma \ref{lem:1} ensure that Algorithm \ref{Alg:unboundedf} is well-defined 
and works properly. More precisely, we have the following result.
\begin{proposition}
Assume Assumption \ref{ass:Lip}. For any choice of $T>0$, $x_0\in \R^{\kappa}$, $ L \geq 1$ and $n_i \in \N$, $c_i \in \R$, $\alpha_i >0$ for $1\leq i\leq L,$ Algorithm \ref{Alg:unboundedf} terminates almost surely within finite time. If additionally $x_0$ is given as in Theorem \ref{thm:1}, the output of Algorithm \ref{Alg:unboundedf} follows the distribution of a Hawkes process with intensity \eqref{eq:intensity}. 
\end{proposition}
It is worth noting that Algorithm \ref{Alg:unboundedf} does not require the sub-criticality condition \eqref{cond:Bremaud}
for non-linear Hawkes processes. Indeed, Algorithm \ref{Alg:unboundedf} applies for instance for the choice $L=1$, $n_1=0$, $c_1=\alpha_1$ and $f(x)=(\mu+x)1_{[0,\infty)}(x)$ with $\mu\geq 0$ for which \eqref{cond:Bremaud} does not hold. 
The only restriction we have to impose is to work with memory kernels which are sum of Erlang kernels, which is a generalization of the approach proposed in \cite{dz}.
In the next section some numerical examples are presented both for bounded and unbounded Lipschitz jump rates $f$. 

\section{Numerical Examples}
\label{sec:numerical_examples}

In this section four numerical examples are given. 
Specifically, we generate first a sample of the Markovian cascade with $L=1$, for a time window $T=100$, order delay $n_1=3$, jump height $c_1=1$, decay rate $\alpha_1=1$ and jump rate $f(x)=(\mu+x)1_{[0,\infty)}(x)$ with $\mu=1$.
We also simulate a Markovian cascade with random jump heights  $c_1$ following a Normal distribution ${\cal N}(0,100)$ and $f(x)=(1+( x/2)^{3/2})\wedge 30$, keeping all others parameters as in the preceding example.
The extension of Algorithm \ref{Alg:unboundedf} for random jump heights is straightforward. 
 Next, we simulate jointly three Markovian cascades with $L=1$
possessing rates of decay $\alpha_1= 0.8, \alpha_1 = 1$ and $\alpha_1 = 1.4$ respectively; in this example $T=500,n_1=3$, the jump heights follow a Normal distribution ${\cal N}(0,100)$ and $f(x)=1+\sigma/(1+e^{-\beta(x-\rho)})$ where $\sigma=20, \beta=1/3$ and $\rho=10$. Finally, we simulate a Markovian cascade for the choices
$L=3$, $n_1=1$, $n_2=3$, $n_3=2$, $\alpha_1=1.3$,
$\alpha_2=0.8$, $\alpha_3=1$, $T=30$, $f(x)=(2+\exp(x/10))\wedge 20$, random jump heights $c_1=c_2=c_3$ following a Normal distribution ${\cal N}(0,25)$.

\begin{figure}[h!]
\includegraphics[scale=0.65]
{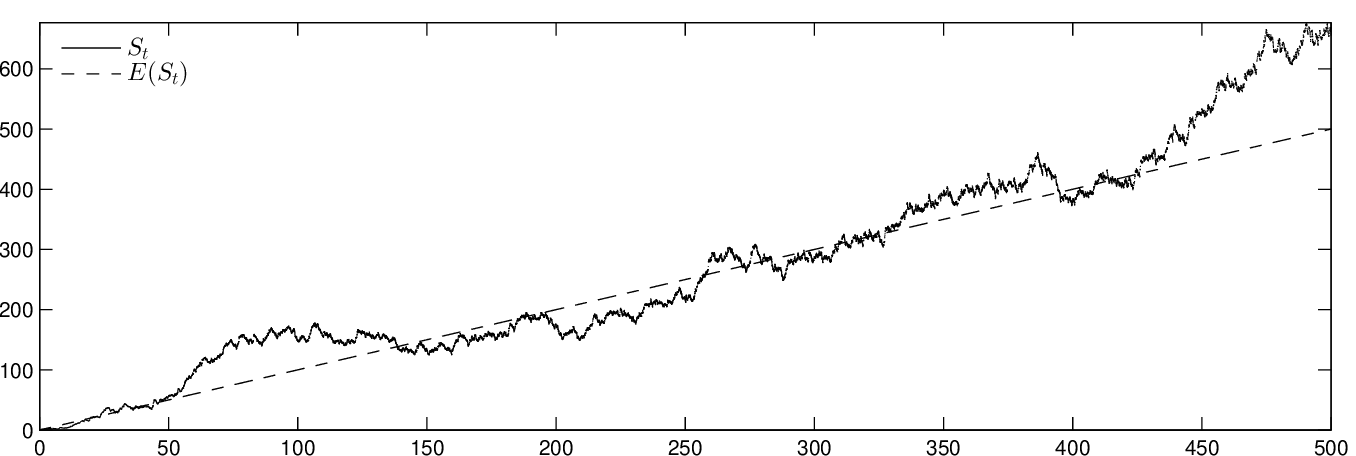}
\hspace{-0.4cm}
\caption{A realization of the process $S=(S_t)_{0\leq t\leq T}$ where $S_t=\sum_{k=0}^{n_1}X^{(k)}_t$  for the choices $L=1$, $n_1=3$, $c_1=\alpha_1=1$, $T=100$ and $f(x)=(1+x/5)1_{[0,\infty)}(x)$ with initial configuration $S_0=0$. The dashed line corresponds to the theoretical mean $E[S_t]=t$ of $S_t$, conditionally on $S_0=0$, obtained by Proposition \ref{prop:mean_sum_X0_Xn}.}
\label{fig:sum_X0_Xn}
\end{figure}

\begin{figure}[h!]
\includegraphics[scale=0.65]{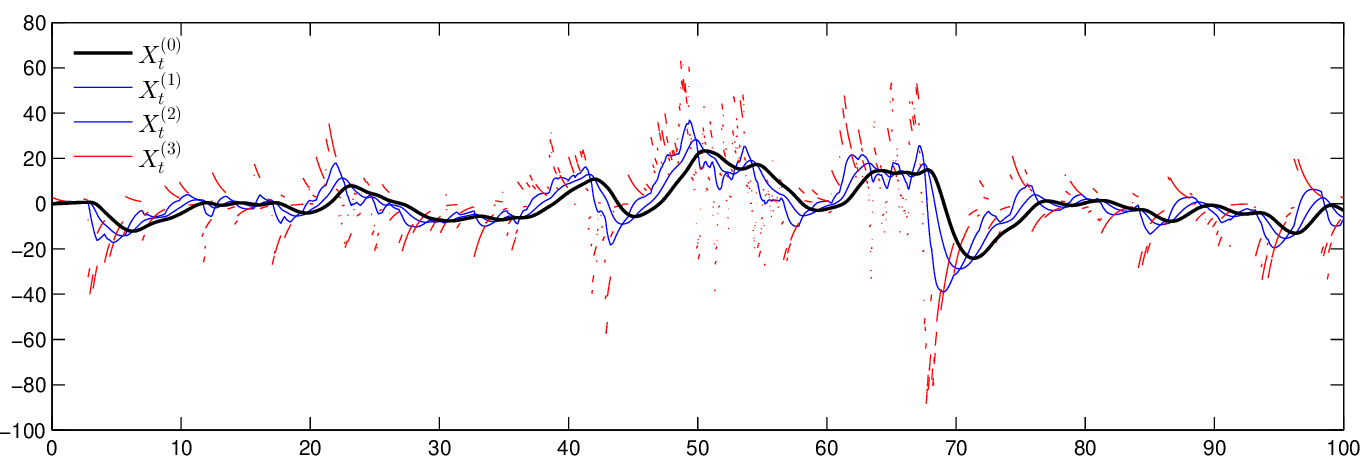}
\hspace{-0.4cm}
\caption{
A realization of the process $X=(X_t)_{0\leq t\leq T}$ with random jump heights following a Normal distribution ${\cal N}(0,100)$ for the choices $L=1, n_1=3$, $\alpha_1=1$, $T=100$, $f(x)=(1+( x/2)^{3/2})\wedge 30$ and $X_0=(0,0,0,0)$. The black trajectory (resp. red) corresponds to the realization of the process $(X_t^{(0)})_{0\leq t\leq T}$ (resp. $(X_t^{(3)})_{0\leq t\leq T}$). The blue trajectories correspond to the realization of the processes $(X_t^{(1)})_{0\leq t\leq T}$
and $(X_t^{(2)})_{0\leq t\leq T}$.}
\label{fig:sum_X0_Xn_with_Normal_distribution}
\end{figure}

\begin{figure}[!h]
\includegraphics[scale=0.65]{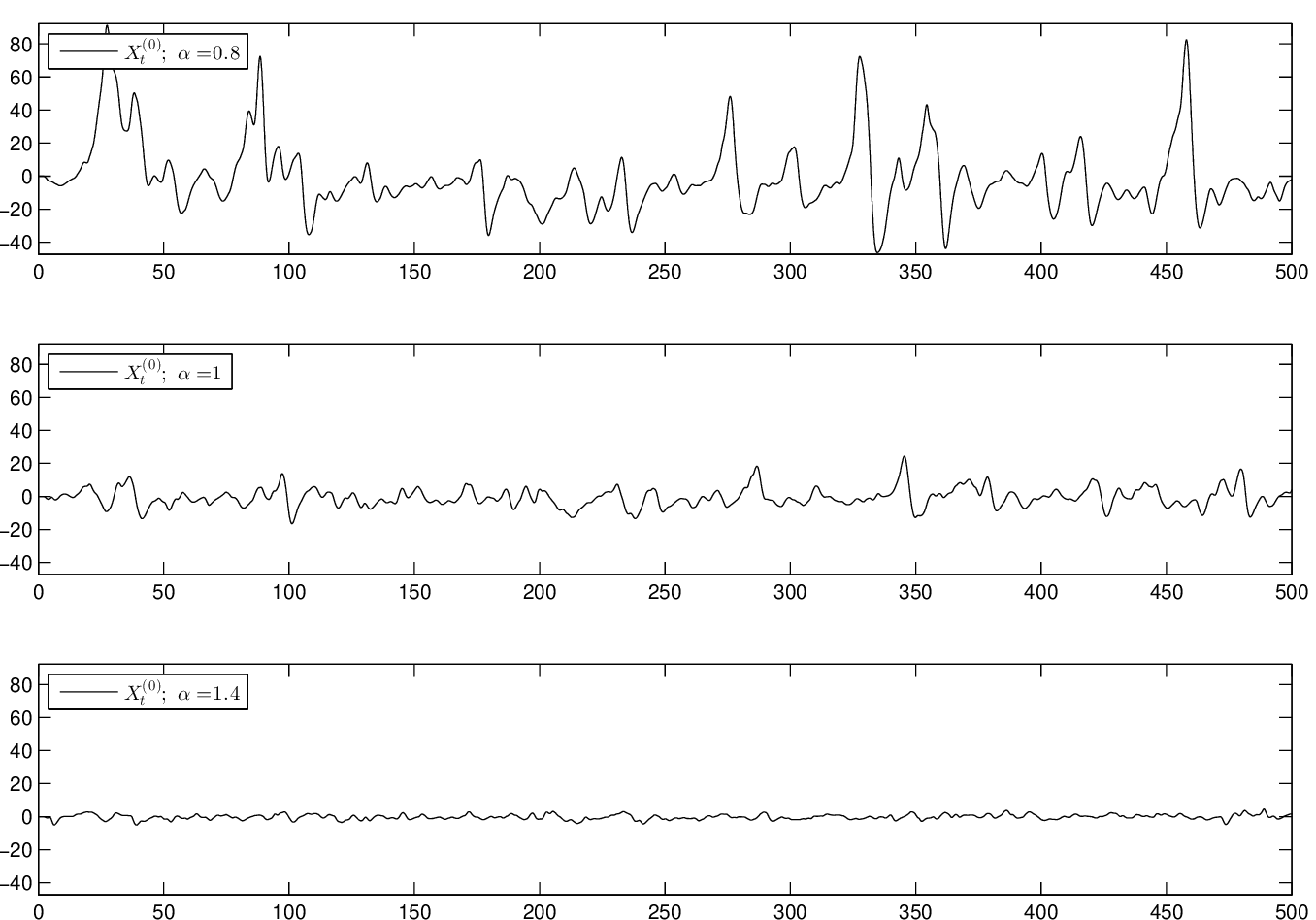} 
\hspace{-0.4cm}
\caption{A joint realization of three Markovian cascades with $\alpha_1=0.8$ (upper panel), $\alpha_1=1$ (middle panel) and $\alpha_1=1.2$ (lower panel) respectively. Here, $L=1$, $T=500,n=3$, the jump heights follow a Normal distribution ${\cal N}(0,100)$ and $f(x)=1+\sigma/(1+e^{-\beta(x-\rho)})$ where $\sigma=20, \beta=1/3$ and $\rho=10$. Notice that the smaller the rate of decay $\alpha$ the larger the oscillations of the process $(X_t^{(0)})_{0\leq t\geq T}$ are.}
\label{fig:three_markovian_cascade}
\end{figure}

\begin{figure}[h!]
\includegraphics[scale=0.65]{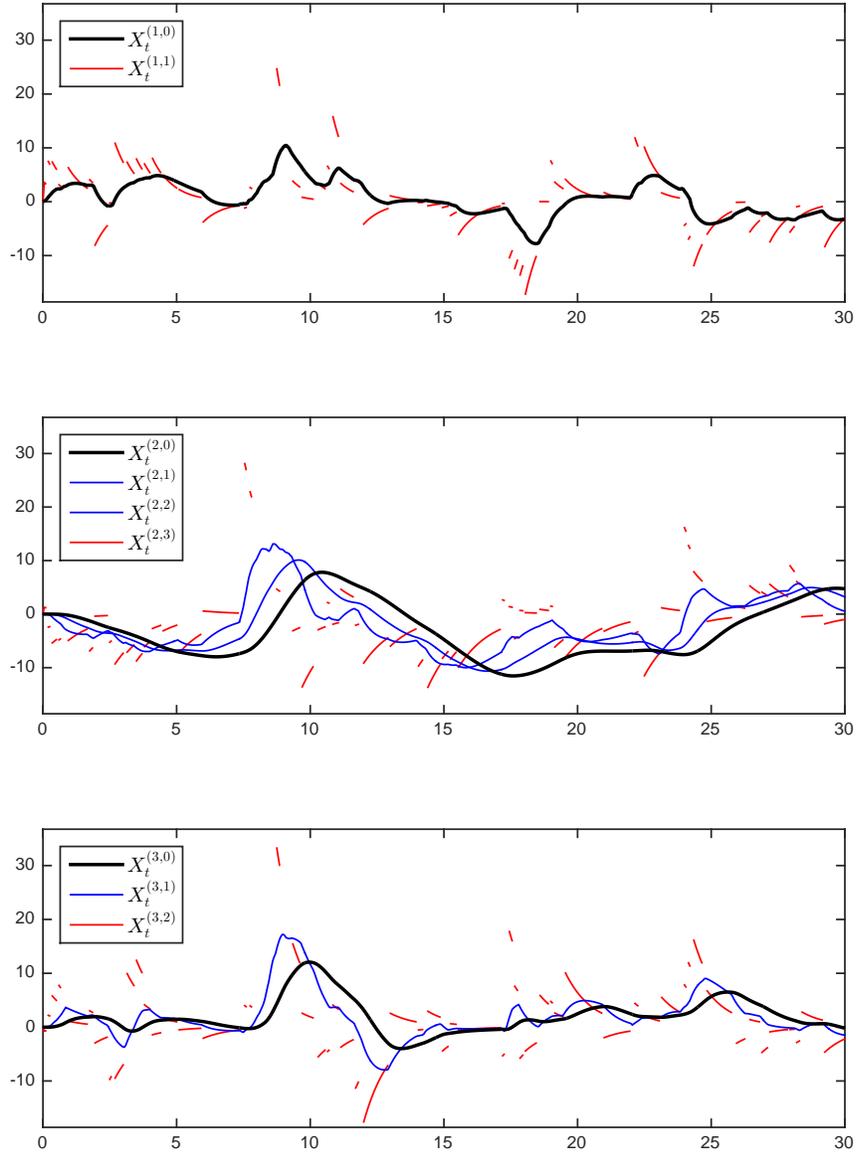}
\caption{A realization of the process $X=(X_t)_{0\leq t\leq T}$ for $L=3$, $n_1=1$, $n_2=3$, $n_3=2$, $\alpha_1=1.3$,
$\alpha_2=0.8$, $\alpha_3=1$, $T=30$, $f(x)=(2+\exp(x/10))\wedge 20$, random jump heights $c_1=c_2=c_3$ following a Normal distribution ${\cal N}(0,25)$ and initial configuration $X_0=0^9$. In the upper panel it is shown separately the realization of the process $(X^{(1)}_t)_{0\leq t\leq T}$, in the middle panel (resp. lower panel) that of $(X^{(2)}_t)_{0\leq t\leq T}$ (resp. $(X^{(3)}_t)_{0\leq t\leq T}$). }
\label{fig:Markovian_cascade_L_3}
\end{figure}

The results are presented in Figures \ref{fig:sum_X0_Xn}, \ref{fig:sum_X0_Xn_with_Normal_distribution}, \ref{fig:three_markovian_cascade} and \ref{fig:Markovian_cascade_L_3} respectively.
In order to test if Algorithm \ref{Alg:unboundedf} works properly, we use the following result.
\begin{proposition}
\label{prop:mean_sum_X0_Xn}
Let $X$ be the Markov process whose generator is given by \eqref{eq:gen} with $L=1, c_1=1, \alpha_1>0, n_1\in\N$   and $f(x)=(\mu+x)1_{[0,\infty)}(x)$ with $\mu\geq 0.$ For any $t\geq 0,$ we 
write $S_t=\sum_{k=0}^{n_1}X^{(k)}_t$. Then for any $x=(x^{(0)},\ldots ,x^{(n_1)})\in \R_+^{n_1  +1 },$ 
\begin{equation}
E[S_t]=\sum_{k=0}^{n_1}x^{(k)}+
\left\lbrace
  \begin{tabular}{cc}
  $\dfrac{\mu}{1-\alpha}\big(e^{t(1-\alpha )}-1\big)$, & \mbox{if} \ $\alpha  \neq 1$  \\
  $\mu t$ & \mbox{if} \ $\alpha =1.$ 
    \end{tabular}
\right.
\end{equation}
\end{proposition}
\begin{proof}
For $\R_+^{n_1 +1}\ni (y^{(0)},\ldots, y^{(n_1)})\mapsto g(y^{(0)},\ldots, y^{(n_1)})=y^{(0)}+\ldots+y^{(n_1)}$ one checks that
$$
Lg(y^{(0)},\ldots, y^{(n_1)})=\mu+(1-\alpha )g(y^{(0)},\ldots, y^{(n_1)}).
$$
By Dynkin's formula it follows that for each $t\geq 0$,
\begin{eqnarray*}
E[S_t]=E[g(X_t)]&=&E[g(X_0)]+\int_{0}^{t}E[Lg(X_s)]ds\\
&=& \sum_{k=0}^{n_1} x^{(k)}+\mu t+(1-\alpha  )\int_{0}^{t}E[S_s]ds,
\end{eqnarray*}
from which it is easy to deduce the result by applying Gronwall's inequality.
\end{proof}

\begin{figure}[h!]
\includegraphics[scale=0.63]{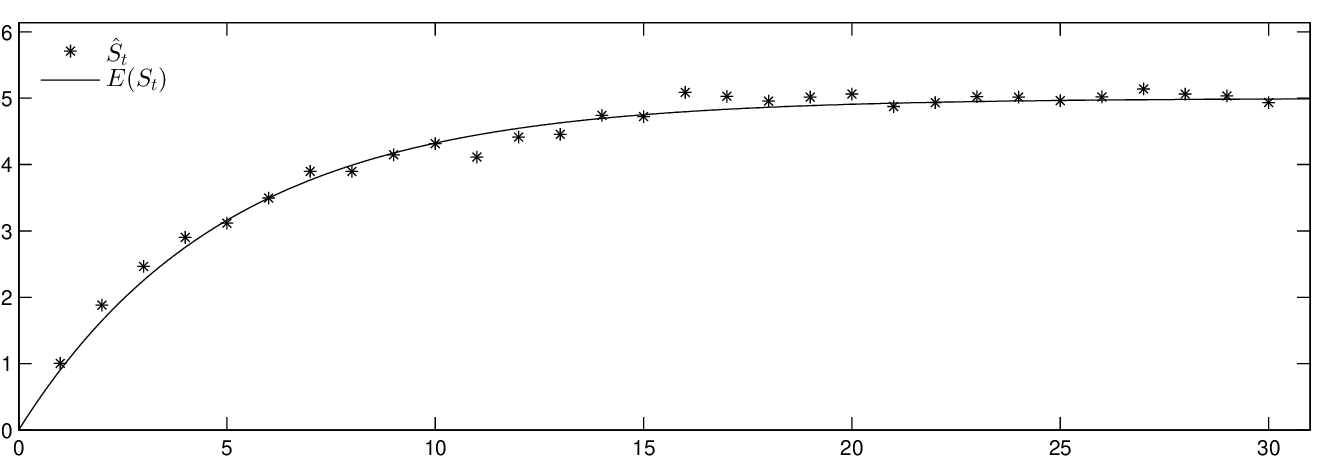}
\caption{The graph of $(E[S_t])_{0\leq t\leq T}$, conditionally on $S_0=0$, for the choices $L=1, n_1=3,c_1=1,\alpha_1=1.2,T=30$ and $f(x)=(1+x/5)1_{[0,\infty)}(x)$. 
The marks $*$ corresponds to the empirical expected value $\hat{S}_t$  of $S_t$ computed at times $t\in \{0,1,\ldots, 30\}$ based on 100 simulated samples of $(S_t)_{0\leq t\leq T}$.}
\label{fig:estimated_sum}
\end{figure}

A comparison between the formula for $E[S_t]$ and and its estimated counterpart denoted by $\hat{S}_t$ is presented in Figure \ref{fig:estimated_sum}.

\section*{Acknowledgements}
The authors thank an anonymous referee for critical reading and useful remarks. E.L. thanks B. Cloez for fruitful discussions concerning Wasserstein coupling at an early stage of this work. 
This research has been conducted as part of the project Labex MME-DII (ANR11-LBX-0023-01); it is part of USP project {\em Mathematics, computation, language
and the brain}, FAPESP project {\em Research, Innovation and
Dissemination Center for Neuromathematics} (grant 2013/07699-0), CNPq projects {\em Stochastic modelling of the brain activity}
(grant 480108/2012-9) and {\em Plasticity in the brain after a brachial plexus lesion} (grant 478537/2012-3).
AD and GO are fully supported by a FAPESP fellowship (grants  2016/17791-9 and 2016/17789-4 respectively). 

\appendix

\section{Proof of Theorem 3 for general $L$}
\label{Appendix}
\begin{proof}
We now prove Theorem \ref{thm:Doeblin} in the case $L> 1.$ 
We write $ X_t = (X_t^{ (1) } , \ldots , X_t^{ (L)}) $ and $ \varphi_t (x) = (\varphi_t^{(1)} ( x), \ldots , \varphi_t^{(L)}  (x) )$ for the flow given in \eqref{eq:flow}. 
Recall that elements of $\R^\kappa $ are denoted by $ x = (x^{(1)}, \ldots, x^{(L)}), $ where $ x^{ (i ) } = (x^{(i, 0)}, \ldots , x^{(i, n_i ) } )$ for each $1\leq i\leq L.$ 

We prove by induction that for all $ 1 \le k \le L, $ for all $x^* \in \R^\kappa $ there exists a neighbourhood $C$ of $x^* $ such that for all $ x \in C,$ 
\begin{multline}\label{eq:ind}
 P_{k T } (x, dy ) \geq \beta_k \nu_k ( y^{(1)}, \ldots , y^{(k)} ) \\
 Q_k ( (x,y^{(1)}, \ldots , y^{(k)} ) , dy^{(k+1)}  \ldots d y^{(L) } ) dy^{(1)} \ldots dy^{(k)} ,
\end{multline} 
where for open sets $ I_l \subset \R^{n_l + 1} , 1 \le l \le k , $ $ \nu_k  $ is the  uniform density on $ I_1 \times \ldots \times I_k, $  and where $Q_k $ is a transition kernel from $ \R^{\kappa } \times \R^{n_1 + \ldots + n_k +k }\to \R^{n_{k+1} + \ldots + n_L + L- k  } .$

{\bf Proof of \eqref{eq:ind} for $k=1.$ }
We proceed as in part I and use the jump noise produced by $n_1 + 1 $ jumps  occurring during $ [0, T ] $ to produce a density for $ X^{(1) }:$ We impose inter-jump waiting times $ t_1, \ldots , t_{n_1 +1} $ under the constraint that $t_1 + \ldots + t_{n_1 + 1 } < T .$ To each jump time we associate jump heights $ c_1, \ldots , c_{n_1 + 1 }  ,$ where each $c_l $ is an element of $ \R^L, $ that is, $c_l = (c_{l}^{( i)}, 1 \le i \le L ) , 1 \le l \le n_1 + 1.$ 

In what follows we shall write  ${\bf \underline{c}}=(c_1, \ldots, c_{n_1+1}) $ and $ {\bf \underline c}^{(i)} = (c_{1 }^{(i)}, \ldots, c_{n_1 + 1 }^{(i)}  ) .$ Moreover, we define  ${\bf\underline{s}}=(s_1, \ldots , s_{n_1+1}),  $ for $ s_k = T - {t_1 + \ldots + t_k } .$ We call $ {\bf \underline s } $ {\it admissible} if $ T > s_1 > \ldots > s_{n_1 + 1 } > 0 .$  

Then, conditionally on $X_0 = x $ and on the above choices,  the position of $ X_T $ is given by 
\begin{equation}\label{def:flowofZbis}
\gamma (x,{\bf \underline{c}} , {\bf \underline{s}}):=\varphi_T  ( x  )+\sum_{i=1}^L \left( c_{1}^{(i)}  e^{- \alpha_i s_1}v_i ({s_1})+\ldots + c_{n_1+1}^{(i)} e^{- \alpha_i s_{n_1 +1}}v_i({s_{n_1+1}}) \right) .
\end{equation}
Here, $v_i ( s) \in \R^\kappa , 1 \le i \le L, $ is the vector given by $ (v_i ( s) )^{(i, k ) } = \frac{ s^{n_i - k }}{(n_i - k )! }, $ for $ 0 \le k \le n_i, $ and with zero entries else. We shall write shortly $ \gamma^{(i)} ({x,{\bf \underline{c}}} , {\bf \underline{s}}) \in \R^{n_i + 1 } $ for the $i-$the coordinate of $ \gamma (x,{\bf \underline{c}}, {\bf \underline{s}}) , $ that is, 
$$\gamma (x,{\bf \underline{c}} , {\bf \underline{s}}) = ( \gamma^{(i)} (x, {\bf \underline{c}}, {\bf \underline{s}}) , 1 \le i \le L ).$$ 

In what follows,  we will use the product form of the flow \eqref{eq:flow}. By this we mean the fact that by the explicit form of the flow in \eqref{eq:flow}, for each $ 1 \le i \le L ,$ we have that 
$$ \varphi^{(i, k)}_t(x )=e^{-\alpha_i t}\sum_{m=0}^{n_i-k}\frac{t^m}{m!}x^{(i, k+m)} =: \varphi^{(i, k)}_t(x^{(i)} ) ;$$
that is $ \varphi^{(i)}_t (x) = \varphi^{(i)}_t (x^{(i)}) $ does only depend on $ x^{(i)} ,$ for any $ 1 \le i \le L.$ As a consequence, 
$$ \gamma^{(i)} (x,{\bf \underline{c}} , {\bf \underline{s}}) =\gamma^{(i)} (x^{(i)},{\bf \underline{c}}^{(i)} , {\bf \underline{s}}) $$
does also depend only on $x^{(i)}$ and on ${\bf \underline c}^{(i)}.$  As usual, we shall write, for any fixed pair $ (x^{(i)},{\bf \underline{c}}^{(i)}),$  
$$ \gamma^{(i)}_{ (x^{(i)},{\bf \underline{c}}^{(i)})} :    {\bf \underline{s}} \mapsto  \gamma^{(i)} (x^{(i)},{\bf \underline{c}}^{(i)} , {\bf \underline{s}}) $$
and similarly, 
$$ \gamma^{(i)}_{ ({\bf \underline{c}^{(i)}} , {\bf \underline{s}})} : x^{(i) } \mapsto \gamma^{(i)} (x^{(i)},{\bf \underline{c}}^{(i)} , {\bf \underline{s}}) .$$

Fix any $x^* \in \R^\kappa $ and fix $ {\bf \underline c^*} $ such that $ {\bf \underline c}_{l}^{(*, i)   } = (a_i + b_i)/2 ,$ for all $ 1 \le l \le n_1 + 1 , $ where $ (a_i, b_i ) \subset \R $ such that $ 0 \notin (a_i, b_i ) $ and $G_i ( (a_i, b_i) ) > 0, $ for all $ 1 \le i \le L.$ 
Then there exists an open neighbourhood $ J_1 = {B}_{R}(x^{(*, 1)} )\times  {B}_R ({\bf \underline{c}}^{(*, 1)} ) $ of the pair $ (x^{(*, 1)}  , {\bf \underline c}^{(*, 1)} ) $ with $ R < \min_i (b_i - a_i  ) , $ and an open set $ I \subset \R^{n_1 + 1 } $ and for any pair $(x^{(1) }, {\bf \underline c}^{(1)}  ) \in  J_{1} $ an open set $ W_{x^{(1)} , {\bf \underline c}^{(1)} } \subset \R_+^{n_1 + 1 } $ such that 
$$
\tilde{\gamma}^{(1)}_{x^{(1)} ,{\bf \underline{c}}^{(1)} }({\bf \underline{s}}):\left\lbrace
\begin{array}{c}
W_{x^{(1)}, {\bf \underline{c}}^{(1)}}\to I \\
{\bf \underline{s}} \mapsto \gamma^{(1)} (x^{(1)} ,{\bf \underline{c}}^{(1)}, {\bf \underline{s}}),
\end{array}
\right.
$$
is a diffeomorphism. Moreover, 
$$ \tilde{\beta_1}= \inf_{ x \in B_R ( x^* ) ,  {\bf \underline c} \in B_R ( {\bf \underline c}^* )   } \inf_{{\bf \underline{s}}\in W_{x^{(1)} ,{\bf \underline{c}}^{(1)} }} q_{x,{\bf \underline{c}}}({\bf \underline{s}}) \Big|\mbox{det}\Big(\frac{\partial \tilde \gamma^{(1)}_{x^{(1)} , {\bf \underline{c}^{(1)} }}({\bf \underline{s}})}{\partial {\bf \underline{s}}}\Big)^{-1}\Big|  > 0 . $$

Let now $ A_i \in {\cal B} ( \R^{n_i + 1 } ) , $ for all $ 1 \le i \le L .$ Then for all $ x \in  {B}_{R }(x^*),$ using the change of variables $y^{(1)} = \tilde \gamma^{(1)}_{ x^{(1)}, {\bf \underline c}^{(1)} } ( {\bf \underline s } )   ,$ 
\begin{multline*}
P_T ( x, A_1 \times \ldots \times A_L ) \geq \tilde{\beta_1} \int_I 1_{A_1} ( y^{(1)}) d y^{(1)} \Big[ \int_{ B_{R}  ( {\bf \underline c}^* )} \prod_{i=1}^L \underline G_i (d {\bf \underline c}^{(i)})   \\
1_{ A_2 } ( \gamma_{ x^{(2)} , {\bf \underline c}^{(2)} }^{(2)}  \circ  (\tilde{\gamma}^{(1)}_{x^{(1)} ,{\bf \underline{c}}^{(1)} })^{-1} (y^{(1)} ) ) \ldots  1_{ A_L \ } ( \gamma_{ x^{(L)} , {\bf \underline c}^{(L)} }^{(L)}  \circ  (\tilde{\gamma}^{(1)}_{x^{(1)} ,{\bf \underline{c}}^{(1)} })^{-1} (y^{(1)} )) \Big ] \\
=: \beta_1 \int_{\R^{n_1 +1} } \nu_1 ( y^{(1) } ) 1_{A_1} ( y^{(1)}) d y^{(1)} \int_{A_2 \times \ldots \times A_L  } Q_1 ((x,  y^{ (1) } ), dy^{(2)} \ldots d y^{(L) } ) ,
\end{multline*}
where $ \nu^{(1)} $ is the uniform density on $I$, $\beta_1=\tilde{\beta}\nu_1(I)$
and  $\underline G_i (d {\bf \underline c}^{(i)})  = \prod_{ l=1}^{n_1 + 1 } G_i ( d c_l^{(i)} ) .$ This proves \eqref{eq:ind} for $k=1.$  

{\bf Induction step : $ k - 1 $ implies $k.$ }
Suppose that  we have already established the result for $k - 1.$  Let $ A_i  \in {\cal B}  ( \R^{n_i +1 } ),$ for all $ 1 \le i \le L.$  We have 
\begin{multline}\label{eq:tobecitedzero}
 P_{kT } ( x, A_1 \times \ldots \times A_L ) =  \int P_T ( y,  A_1 \times \ldots \times A_L ) P_{(k-1) T}  (x, dy )  \, \\
 \geq \beta_{k-1}  1_{C } ( x) 
 \int   P_T ( y,  A_1 \times \ldots \times A_L ) \nu_{k-1} (y^{(1 )}, \ldots , y^{(k-1)} )\\
  Q_{k-1}  ( (x, y^{(1)} , \ldots , y^{(k-1)} ), d y^{(k)} \ldots dy^{(L) } )  dy^{(1)} \ldots dy^{(k-1)} .
\end{multline}
We work conditionally on the choice of $ y = (y^{(1)}, \ldots , y^{(L)} ) $ and proceed as in the first part, using the jump noise (of a sufficient number of jumps) to create a density for the variable $y^{(k)} $ and proving that the already produced density $ \nu_{k-1} (y^{(1 )}, \ldots , y^{(k-1)} )$ of the first $k-1$ variables is well preserved.  

As in the first step, we start with 
\begin{equation}\label{eq:tobecitedfirst}
 P_T ( y, A_1 \times  \ldots \times A_L ) \geq E_y \left( \prod_{i=1}^L 1_{A_i} ( X^{(i)}_T) , N_T = n_k +1 \right)  .
\end{equation} 
We impose inter-jump waiting times $ t_1, \ldots , t_{n_k +1} $ under the constraint that $t_1 + \ldots + t_{n_k + 1 } < T $ and associated jump heights $ c_1, \ldots , c_{n_k + 1 }  ,$ where each $c_l $ is an element of $ \R^L, $ that is, $c_l = (c_l^{ ( i) }, 1 \le i \le L ) .$ 
These inter-jump waiting times will produce a density for the $k-$th variable $X_T^{(k)} $ as in the preceding steps. 

To start with, let us introduce the following notation.
For all $ 1 \le l < m \le L, $ we write  $ y^{l: m} := ( y^{(l)}, \ldots , y^{(m)}) ,$ $ d  y^{l: m} :=  dy^{(l)} \ldots dy^{(m)},  $${\bf \underline c}^{l: m } := ({\bf \underline c}^{(l)}, \ldots , {\bf \underline c}^{(m)} ) .$  For $x\in C$ we have
\begin{multline}\label{eq:tobecitedbis}
 \int P_T ( y, A_1 \times \ldots \times  A_L )  \nu_{k-1} (y^{ 1: k-1} )  Q_{k-1} ( (x, y^{1: k-1)} ), d y^{k:L} )  dy^{1: k-1} \\ 
\geq \int \nu_{k-1} (y^{ 1: k-1} ) Q_{k-1} ( (x,y^{ 1: k-1} ), d y^{k: L} )  d y^{ 1: k-1}  \int \prod_{i=1}^L \underline G_i (d {\bf \underline c}^{(i)} )  \\
\int ds_{1}\ldots ds_{n_k+1} 1_{A_1\times \ldots \times A_{k-1}}(\gamma^{1: k-1} (y^{1:k-1} , {\bf \underline{c}}^{1:k-1} , {\bf \underline{s}}) ) 1_{A_k}(\gamma^{( k) } (y^{(k )} , {\bf \underline{c}}^{(k)} , {\bf \underline{s}})) \\
 1_{A_{k+1}\times \ldots \times A_{L}}(\gamma^{k+1: L } ( y^{k+1:L} , {\bf \underline{c}}^{k+1:L},  {\bf \underline{s}}))  q_{y,{\bf \underline{c}}}({\bf \underline{s}})  .
\end{multline}
  
Let $ N_1 := \sum_{ l=1}^{k- 1} (n_l +1) $ be the dimension of $ y^{1 : k-1} .$ We introduce now for any fixed $ {\bf \underline c^{  1: k }  } $ having all entries non zero and  $y^{(k)} \in \R^{n_k +1}, $ 
$$ \Phi_{ (    {\bf \underline c^{  1: k }  } , y^{(k)}) }  :\left\lbrace
\begin{array}{lll}
\R^{ N_1 } \times \R^{ n_k +1 } &\to &\R^{ N_1 + n_k + 1 } \\
( y^{ 1: k-1} , {\bf \underline s} )& \mapsto &  (\gamma^{1: k-1} (y^{1:k-1} , {\bf \underline{c}}^{1:k-1} , {\bf \underline{s}}) , \gamma^{( k) } (y^{(k )} , {\bf \underline{c}}^{(k)} , {\bf \underline{s}}) ) . 
\end{array}  \right. $$  
We write 
$$
\frac{\partial \Phi_{ (    {\bf \underline c^{  1: k }  } , y^{(k)}) } ( y^{ 1: k-1} , {\bf \underline s} )}{\partial y^{1: k-1 } \partial {\bf \underline{s}}}=\Big[\frac{\partial \Phi_{ (    {\bf \underline c^{  1: k }  } , y^{(k)}) } ( y^{ 1: k-1} , {\bf \underline s} )}{\partial y^{(1)} },\ldots, \frac{\partial \Phi_{ (    {\bf \underline c^{  1: k }  } , y^{(k)}) } ( y^{ 1: k-1} , {\bf \underline s} )}{\partial s_{n_k+1 } } \Big]
$$ 
to denote the Jacobian matrix of the the map $(y^{1: k-1} , \underline{\bf s}) \mapsto \Phi_{ (    {\bf \underline c^{  1: k }  } , y^{(k)}) } (y^{1: k-1} , \underline{\bf s}). $ By the properties of the flow \eqref{eq:flow}, it follows that 
$$ \frac{\partial \Phi_{ (    {\bf \underline c^{  1: k }  } , y^{(k)}) } ( y^{ 1: k-1} , {\bf \underline s} )}{\partial y^{(1)} } =
\left(  \begin{array}{c}
\frac{ \partial \gamma^{(1)}  ( y^{(1)}, {\bf \underline {c}^{(1)} , {\bf \underline s}  ) }}{\partial y^{(1)} } \\
0 \\
\vdots \\
0 
\end{array} 
\right) 
= 
 \left(  \begin{array}{c}
\frac{ \partial \varphi_T ( y^{(1)} ) }{\partial y^{(1)} } \\
0 \\
\vdots \\
0 
\end{array} 
\right) ,$$
where by the ``cascade structure'' of the flow \eqref{eq:flow}, 
$$ 
\frac{ \partial \varphi_T ( y^{(1)} ) }{\partial y^{(1)} } = 
\left(  \begin{array}{cccc} 
e^{-  \alpha_1 T }   & * &*&* \\
0 & e^{-  \alpha_1 T }  &*&* \\
\vdots & &\ddots  & *\\
0 & \cdots  & 0 & e^{-  \alpha_1 T } 
\end{array} 
\right) .$$
Therefore, 
$$  \frac{\partial \Phi_{ (    {\bf \underline c^{  1: k }  } , y^{(k)}) } ( y^{ 1: k-1} , {\bf \underline s} )}{\partial y^{1: k-1 } \partial {\bf \underline{s}}} = \left( 
\begin{array}{cc}
A & * \\
0 & B 
\end{array} 
\right) , $$
where $A$ is an $N_1 \times N_1 $ upper diagonal matrix having entries of the type $ e^{ - \alpha_i T } , 1 \le i \le k-1 ,$ on the diagonal, where $ 0 $ is the $ (n_k+1) \times N_1 $ matrix having all zero entries, and where 
$B = \left( \frac{ \partial  \gamma^{( k) } (y^{(k )} , {\bf \underline{c}}^{(k)} , {\bf \underline{s}})}{\partial \bf \underline s} \right).$ According to Part I of the proof, $ \det B \neq 0 $ whenever $ {\bf \underline{c}}^{(k)} $ has only non zero coordinates, for any fixed $ y^{(k)}, $ for admissible $ {\bf \underline s }.$ It is immediate to see that 
$$ \det \left(  \frac{\partial \Phi_{ (    {\bf \underline c^{  1: k }  } , y^{(k)}) } ( y^{ 1: k-1} , {\bf \underline s} )}{\partial y^{1: k-1 } \partial {\bf \underline{s}}} \right) = e^{ - ( \alpha_1 (n_1 +1) + \ldots \alpha_{k-1} (n_{k-1} +1 ) ) T } \cdot \det B \neq 0 $$ 
in this case. 
Therefore, arguing as in Part I,  for any $ y^{(*,k) }  \in \R^{n_k +1}   $ and $ {\bf \underline c^{ (*, 1: k)}  } $ having non zero entries, for all $ y^{(*, 1 : k-1)} $ and an admissible $ {\bf \underline {s}^*} $  there exist open neighbourhoods $ J_{k } =  B_{R } (y^{ (*, k) }) \times  B_{R } ({\bf \underline c }^{(*, 1: k )} ) $ and an open set $ I \subset \R^{N_1 + n_k +1} $ containing $ (  y^{(*, 1 : k-1)}  , {\bf \underline {s}^*} )$  and for any pair $ (y^{(k)}, {\bf \underline c}^{(1: k)} ) \in  J_{k} $ an open set $ W_{y^{(k)}, {\bf \underline c}^{(1: k)}} $ such that 
$$
\tilde{\Phi}_{ (    {\bf \underline c^{  1: k }  } , y^{(k)}) } (y^{1:k-1}, {\bf \underline{s}}):\left\lbrace
\begin{array}{ll}
W_{y^{(k)},{\bf \underline{c}^{(1: k)} } } &\to I \\
(y^{1:k-1}, {\bf \underline{s}}) & \mapsto \Phi_{ (    {\bf \underline c^{  1: k }  } , y^{(k)}) } ( y^{ 1: k-1} , {\bf \underline s} ) ,
\end{array}
\right.
$$
is a diffeomorphism. 

In what follows, in order to ease notation, we shall shortly write 
$$ \tilde \Phi (y^{1:k-1}, {\bf \underline{s}}) := \tilde{\Phi}_{ (    {\bf \underline c^{  1: k }  } , y^{(k)}) } (y^{1:k-1}, {\bf \underline{s}}) $$
and 
\begin{equation}\label{eq:inverseouf}
 \Psi =  \tilde \Phi ^{-1} 
\end{equation} 
for the associated inverse function. $\Psi $ taking values in some (subset of) $\R^{N_1 + n_k +1 }, $ we shall write as before $ \Psi^{(i)} $ for its coordinates and $ \Psi^{1: N_1 } $ for the first $N_1 $ of its coordinates, corresponding to  $y^{1: k-1} ,$ and $\Psi^{N_1+1 : N_1+n_k +1 } $ for the last coordinates, corresponding to $ {\bf \underline s}. $  

We choose any  $x\in C$ and  $y^* \in $ supp $(  \nu_{k-1} (y^{1:k-1})   Q_{k-1} ( (x, y^{1: k-1}) , d y^{k:L}  ) d y^{1:k-1} )   $  and obtain (recall  \eqref{eq:tobecitedfirst} and \eqref{eq:tobecitedbis}),
\begin{multline*}
 \int P_T ( y, A_1 \times \ldots \times  A_L )  \nu_{k-1} (y^{ 1: k-1} )  Q_{k-1} ( (x, y^{1: k-1)} ), d y^{k:L} )  dy^{1: k-1} \\ 
\geq \int_{ B_{R  } (y^*)}  \nu_{k-1} (y^{ 1: k-1} ) Q_{k-1} ( (x,y^{ 1: k-1} ), d y^{k: L} )  d y^{ 1: k-1}  \int_{{B}_{R}({\bf \underline{c}}^*)}\prod_{i=1}^L \underline G_i (d {\bf \underline c}^{(i)} )  \\
\int_{\R_+^{n_k +1} }   1_{ W_{y^{(k)},{\bf \underline{c}^{(1: k)} }}}  ( y^{1:k-1}, {\bf \underline s})  \; q_{y,{\bf \underline{c}}}({\bf \underline{s}})\; 1_{A_1\times \ldots \times A_{k} }(\tilde \Phi (y^{1:k-1}, {\bf \underline{s}})  ) \\
  1_{ A_{k+1} \times \ldots A_L} ( \gamma^{k+1: L }_{y^{k+1:L} , {\bf \underline{c}}} ({\bf \underline{s}}))ds_{1}\ldots ds_{n_k+1} .
\end{multline*}
In the above formula, $R$ is chosen sufficiently small such that  ${B}_{R}({\bf \underline{c}}^*)  $ contains only jump heights with non zero entries. We then use the change of variables  
$$ z^{1: k } :  =\tilde \Phi (y^{1:k-1}, {\bf \underline{s}}) . $$ 
Choose now $R$ sufficiently small such that 
\begin{eqnarray*}
\tilde \beta_k &:= & \inf_{y \in B_R (y^*) } \inf_{ {\bf \underline{c}}\in B_R ({\bf \underline{c}}^*) } \;  \inf_{{\bf \underline{s}}\, : \, (y^{1:k-1}, {\bf \underline s})   \in W_{y^{(k)},{\bf \underline{c}^{1: k }}}} q_{y,{\bf \underline{c}}}({\bf \underline{s}}) \times  \\
&&\times  \Big|\mbox{det}\Big(\frac{\partial \tilde \Phi (y^{1:k-1}, {\bf \underline{s}})}{\partial y^{1: k-1 } \partial {\bf \underline{s}}}\Big)^{-1}\Big| \nu_{k-1} (y^{ 1: k-1} )   > 0   .
\end{eqnarray*}
Let $C_2^* =  B_{R  } (y^{ (*, k)}, \ldots , y^{ ( *, L)}) . $  Then
\begin{eqnarray*}
&&\int P_T ( y, A_1 \times \ldots \times A_L  )  \nu_{k-1} (y^{ 1: k-1} )  Q_{k-1} ( (x, y^{1: k-1)} ), d y^{k:L} )  dy^{1: k-1} \\ 
&&\geq \tilde \beta_k \int_{I}  ( \prod_{j=1}^k 1_{ A_j } (z^{ ( j)}  ) d z^{1: k }    \Big[   \int_{{B}_R ({\bf \underline{c}}^*)}\prod_{i=1}^L \underline G_i (d {\bf \underline c}^{(i)} )  \\
&&\int_{C_2^* }  Q_{k-1} (( x, \Psi^{1:N_1} ( z^{1:k}  ) , dy^{k:L} )   1_{A_{k+1} \times \ldots A_L } (  \gamma^{k+1: L }_{y^{k+1:L} , {\bf \underline{c}}}\circ  \Psi^{N-1 + 1 : N_1 + n_k + 1 } )  (z^{1:k} )  
\Big]  \\
&&=: \tilde \beta_k \int_{I}  ( \prod_{j=1}^k 1_{ A_j } (z^{ ( j)}  ) d z^{1: k } \left[  \int Q_k ( (x, z^{1:k} ), d z^{ k+1 : L } ) 1_{ A_{k+1} } ( z^{(k+1 ) } )\cdot \ldots \cdot 1_{A_L} ( z^{(L) } ) \right] . 
\end{eqnarray*}
Together with \eqref{eq:tobecitedzero}, this shows that \eqref{eq:ind} holds also for $ k ,$ and this finishes the induction step.  By taking finally  $ k = L $ in \eqref{eq:ind}, this implies the assertion of the Theorem. 

\end{proof}

\bibliography{Bibli}{}

\begin{thebibliography}{23}

\bibitem[\protect\citeauthoryear{Az{\'e}ma, Duflo and Revuz}{1969}]{adr}
\begin{barticle}[author]
\bauthor{\bsnm{Az{\'e}ma},~\bfnm{J.}\binits{J.}},
  \bauthor{\bsnm{Duflo},~\bfnm{M.}\binits{M.}} \AND
  \bauthor{\bsnm{Revuz},~\bfnm{D.}\binits{D.}}
(\byear{1969}).
\btitle{Mesures invariantes des processus de Markov r{\'e}currents}.
\bjournal{S\'em. Proba. III, Lecture Notes in Math.}
\bvolume{88}
\bpages{24-33}.
\end{barticle}
\endbibitem

\bibitem[\protect\citeauthoryear{Bena{\"i}m et~al.}{2015}]{benaim2015}
\begin{barticle}[author]
\bauthor{\bsnm{Bena{\"i}m},~\bfnm{Michel}\binits{M.}},
  \bauthor{\bsnm{Le~Borgne},~\bfnm{Stéphane}\binits{S.}},
  \bauthor{\bsnm{Malrieu},~\bfnm{Florent}\binits{F.}} \AND
  \bauthor{\bsnm{Zitt},~\bfnm{Pierre-André}\binits{P.-A.}}
(\byear{2015}).
\btitle{Qualitative properties of certain piecewise deterministic Markov
  processes}.
\bjournal{Ann. Inst. H. Poincaré Probab. Statist.}
\bvolume{51}
\bpages{1040--1075}.
\bdoi{10.1214/14-AIHP619}
\end{barticle}
\endbibitem

\bibitem[\protect\citeauthoryear{Bolley}{2008}]{bolley}
\begin{barticle}[author]
\bauthor{\bsnm{Bolley},~\bfnm{F.}\binits{F.}}
(\byear{2008}).
\btitle{Separability and completeness for the Wasserstein distance}.
\bjournal{S\'em. Proba. XLI, Lecture Notes in Math.}
\bvolume{1934}
\bpages{371-377}.
\end{barticle}
\endbibitem

\bibitem[\protect\citeauthoryear{Bonde~Raad, Ditlevsen and
  L\"ocherbach}{2018}]{madsevasusanne}
\begin{barticle}[author]
\bauthor{\bsnm{Bonde~Raad},~\bfnm{M}\binits{M.}},
  \bauthor{\bsnm{Ditlevsen},~\bfnm{S.}\binits{S.}} \AND
  \bauthor{\bsnm{L\"ocherbach},~\bfnm{E.}\binits{E.}}
(\byear{2018}).
\btitle{Age dependent {H}awkes process}.
\bjournal{arXiv:1806.06370}.
\end{barticle}
\endbibitem

\bibitem[\protect\citeauthoryear{Br{\'e}maud and Massouli{\'e}}{1996}]{bm}
\begin{barticle}[author]
\bauthor{\bsnm{Br{\'e}maud},~\bfnm{P.}\binits{P.}} \AND
  \bauthor{\bsnm{Massouli{\'e}},~\bfnm{L.}\binits{L.}}
(\byear{1996}).
\btitle{Stability of nonlinear Hawkes processes}.
\bjournal{The Annals of Probability}
\bvolume{24}
\bpages{1563-1588}.
\end{barticle}
\endbibitem

\bibitem[\protect\citeauthoryear{Chevallier}{2017}]{chevallier}
\begin{barticle}[author]
\bauthor{\bsnm{Chevallier},~\bfnm{J.}\binits{J.}}
(\byear{2017}).
\btitle{Mean-field limit of generalized Hawkes processes}.
\bjournal{Stoch. Proc. Appl.}
\bvolume{127}
\bpages{3870-3912}.
\end{barticle}
\endbibitem

\bibitem[\protect\citeauthoryear{Chevallier et~al.}{2015}]{ccdr}
\begin{barticle}[author]
\bauthor{\bsnm{Chevallier},~\bfnm{J.}\binits{J.}},
  \bauthor{\bsnm{Caceres},~\bfnm{MJ.}\binits{M.}},
  \bauthor{\bsnm{Doumic},~\bfnm{M.}\binits{M.}} \AND
  \bauthor{\bsnm{Reynaud-Bouret},~\bfnm{P.}\binits{P.}}
(\byear{2015}).
\btitle{Microscopic approach of a time elapsed neural model}.
\bjournal{Math. Mod. \& Meth. Appl. Sci.}
\bvolume{25}
\bpages{2669 - 2719}.
\end{barticle}
\endbibitem

\bibitem[\protect\citeauthoryear{Dassios and Zhao}{2013}]{dz}
\begin{barticle}[author]
\bauthor{\bsnm{Dassios},~\bfnm{Angelos}\binits{A.}} \AND
  \bauthor{\bsnm{Zhao},~\bfnm{Hongbiao}\binits{H.}}
(\byear{2013}).
\btitle{Exact simulation of Hawkes process with exponentially decaying
  intensity}.
\bjournal{Electron. Commun. Probab.}
\bvolume{18}
\bpages{13 pp.}
\bdoi{10.1214/ECP.v18-2717}
\end{barticle}
\endbibitem

\bibitem[\protect\citeauthoryear{Delattre, Fournier and Hoffmann}{2016}]{dfh}
\begin{barticle}[author]
\bauthor{\bsnm{Delattre},~\bfnm{S.}\binits{S.}},
  \bauthor{\bsnm{Fournier},~\bfnm{N.}\binits{N.}} \AND
  \bauthor{\bsnm{Hoffmann},~\bfnm{M.}\binits{M.}}
(\byear{2016}).
\btitle{Hawkes processes on large networks}.
\bjournal{Ann. App. Probab.}
\bvolume{26}
\bpages{216 - 261}.
\end{barticle}
\endbibitem

\bibitem[\protect\citeauthoryear{Ditlevsen and L\"ocherbach}{2017}]{SusanneEva}
\begin{barticle}[author]
\bauthor{\bsnm{Ditlevsen},~\bfnm{S.}\binits{S.}} \AND
  \bauthor{\bsnm{L\"ocherbach},~\bfnm{E.}\binits{E.}}
(\byear{2017}).
\btitle{Multi-class oscillating systems of interacting neurons}.
\bjournal{Stoch. Proc. Appl.}
\bvolume{127}
\bpages{1840-1869}.
\end{barticle}
\endbibitem

\bibitem[\protect\citeauthoryear{Ditlevsen, Yip and
  Holstein-Rathlou}{2005}]{Ditlevsen2005}
\begin{barticle}[author]
\bauthor{\bsnm{Ditlevsen},~\bfnm{S.}\binits{S.}},
  \bauthor{\bsnm{Yip},~\bfnm{K.~P.}\binits{K.~P.}} \AND
  \bauthor{\bsnm{Holstein-Rathlou},~\bfnm{N.~H}\binits{N.~H.}}
(\byear{2005}).
\btitle{Parameter estimation in a stochastic model of the tubuloglomerular
  feedback mechanism in a rat nephron}.
\bjournal{Math. Biosci.\/.}
\bvolume{194}
\bpages{49-69}.
\end{barticle}
\endbibitem

\bibitem[\protect\citeauthoryear{Down, Meyn and Tweedie}{1995}]{dmt}
\begin{barticle}[author]
\bauthor{\bsnm{Down},~\bfnm{D.}\binits{D.}},
  \bauthor{\bsnm{Meyn},~\bfnm{S.~P.}\binits{S.~P.}} \AND
  \bauthor{\bsnm{Tweedie},~\bfnm{R.~L.}\binits{R.~L.}}
(\byear{1995}).
\btitle{Exponential and Uniform Ergodicity of Markov Processes}.
\bjournal{Ann. Probab.}
\bvolume{23}
\bpages{1671--1691}.
\bdoi{10.1214/aop/1176987798}
\end{barticle}
\endbibitem

\bibitem[\protect\citeauthoryear{Hansen, Reynaud-Bouret and
  Rivoirard}{2015}]{hrbr}
\begin{barticle}[author]
\bauthor{\bsnm{Hansen},~\bfnm{N.}\binits{N.}},
  \bauthor{\bsnm{Reynaud-Bouret},~\bfnm{P.}\binits{P.}} \AND
  \bauthor{\bsnm{Rivoirard},~\bfnm{V.}\binits{V.}}
(\byear{2015}).
\btitle{Lasso and probabilistic inequalities for multivariate point processes}.
\bjournal{Bernoulli}
\bvolume{21}
\bpages{83 - 143}.
\end{barticle}
\endbibitem

\bibitem[\protect\citeauthoryear{Hawkes}{1971}]{Hawkes1971}
\begin{barticle}[author]
\bauthor{\bsnm{Hawkes},~\bfnm{Alan~G.}\binits{A.~G.}}
(\byear{1971}).
\btitle{{Spectra of Some Self-Exciting and Mutually Exciting Point Processes}}.
\bjournal{Biometrika}
\bvolume{58}
\bpages{83--90}.
\end{barticle}
\endbibitem

\bibitem[\protect\citeauthoryear{Hawkes and Oakes}{1974}]{ho}
\begin{barticle}[author]
\bauthor{\bsnm{Hawkes},~\bfnm{A.~G.}\binits{A.~G.}} \AND
  \bauthor{\bsnm{Oakes},~\bfnm{D.}\binits{D.}}
(\byear{1974}).
\btitle{A cluster process representation of a self-exciting process}.
\bjournal{J. Appl. Probab.}
\bvolume{11}
\bpages{93 - 503}.
\end{barticle}
\endbibitem

\bibitem[\protect\citeauthoryear{Jacod}{1975}]{Jacod75}
\begin{barticle}[author]
\bauthor{\bsnm{Jacod},~\bfnm{J.}\binits{J.}}
(\byear{1975}).
\btitle{Multivariate {P}oint {P}rocesses: {P}redictable {P}rojection,
  {R}adon-{N}ikodym {D}erivatives, {R}epresentation of {M}artingales}.
\bjournal{Z. Wahrscheinlichkeitstheorie verw. Gebiete}
\bvolume{31}
\bpages{235-253}.
\end{barticle}
\endbibitem

\bibitem[\protect\citeauthoryear{Kammler}{1976}]{kammler}
\begin{barticle}[author]
\bauthor{\bsnm{Kammler},~\bfnm{D.~W.}\binits{D.~W.}}
(\byear{1976}).
\btitle{Approximation with sums of exponentials in $L_p [0, \infty ) $}.
\bjournal{J. of Approx. Theory}
\bvolume{16}.
\end{barticle}
\endbibitem

\bibitem[\protect\citeauthoryear{L\"ocherbach}{2017}]{evalast}
\begin{barticle}[author]
\bauthor{\bsnm{L\"ocherbach},~\bfnm{E.}\binits{E.}}
(\byear{2017}).
\btitle{Convergence to equilibrium for time inhomogeneous jump diffusions with
  state dependent jump intensity}.
\bjournal{arXiv:1712.03507}.
\end{barticle}
\endbibitem

\bibitem[\protect\citeauthoryear{L{\"o}cherbach and
  Loukianova}{2008}]{dashaeva}
\begin{barticle}[author]
\bauthor{\bsnm{L{\"o}cherbach},~\bfnm{Eva}\binits{E.}} \AND
  \bauthor{\bsnm{Loukianova},~\bfnm{Dasha}\binits{D.}}
(\byear{2008}).
\btitle{On Nummelin splitting for continuous time Harris recurrent Markov
  processes and application to kernel estimation for multi-dimensional
  diffusions}.
\bjournal{Stochastic Processes and their Applications}
\bvolume{118}
\bpages{1301 - 1321}.
\end{barticle}
\endbibitem

\bibitem[\protect\citeauthoryear{Meyn and Tweedie}{1993}]{Meyn93}
\begin{barticle}[author]
\bauthor{\bsnm{Meyn},~\bfnm{S.~P.}\binits{S.~P.}} \AND
  \bauthor{\bsnm{Tweedie},~\bfnm{R.~L.}\binits{R.~L.}}
(\byear{1993}).
\btitle{Stability of {Markovian processes III : Foster-Lyapunov} criteria for
  continuous-time processes.}
\bjournal{Adv. Appl. Probab.}
\bvolume{25}
\bpages{487-548}.
\end{barticle}
\endbibitem

\bibitem[\protect\citeauthoryear{Rachev}{1991}]{Rachev}
\begin{bbook}[author]
\bauthor{\bsnm{Rachev},~\bfnm{S.~T.}\binits{S.~T.}}
(\byear{1991}).
\btitle{Probability metrics and the stability of stochastic models}.
\bpublisher{John Wiley and Sons}, \baddress{Chichester, USA}.
\end{bbook}
\endbibitem

\bibitem[\protect\citeauthoryear{Skeldon and Purvey}{2005}]{KidneyDelay}
\begin{barticle}[author]
\bauthor{\bsnm{Skeldon},~\bfnm{A.~C.}\binits{A.~C.}} \AND
  \bauthor{\bsnm{Purvey},~\bfnm{I.~.}\binits{I.~.}}
(\byear{2005}).
\btitle{The effect of different forms for the delay in a model of the nephron}.
\bjournal{Math. Biosci. Eng.}
\bvolume{2(1)}
\bpages{97-109}.
\end{barticle}
\endbibitem

\bibitem[\protect\citeauthoryear{Zhu}{2015}]{zhu2015}
\begin{barticle}[author]
\bauthor{\bsnm{Zhu},~\bfnm{Lingjiong}\binits{L.}}
(\byear{2015}).
\btitle{Large deviations for Markovian nonlinear Hawkes processes}.
\bjournal{Ann. Appl. Probab.}
\bvolume{25}
\bpages{548--581}.
\bdoi{10.1214/14-AAP1003}
\end{barticle}
\endbibitem

\end{thebibliography}
\bibliographystyle{imsart-nameyear}
\end{document}